\pdfoutput=1
\documentclass[11pt]{article}

\usepackage[a4paper,margin=32mm]{geometry}

\usepackage[T1]{fontenc}
\usepackage[utf8]{inputenc}
\usepackage[british]{babel}
\usepackage{lmodern}
\usepackage{microtype}
\usepackage{amsmath,amssymb,amsthm,mathtools}
\numberwithin{equation}{section}
\usepackage{aliascnt}
\usepackage{mathrsfs}

\usepackage{booktabs}
\usepackage{array}
\usepackage{enumitem}

\usepackage{cite}
\usepackage[colorlinks=true,linkcolor=blue,citecolor=blue,urlcolor=blue]{hyperref}
\usepackage[nameinlink,capitalize]{cleveref}
\theoremstyle{plain}
\newtheorem{theorem}{Theorem}[section]
\newaliascnt{lemma}{theorem}
\newtheorem{lemma}[lemma]{Lemma}
\aliascntresetthe{lemma}
\newaliascnt{proposition}{theorem}
\newtheorem{proposition}[proposition]{Proposition}
\aliascntresetthe{proposition}
\newaliascnt{corollary}{theorem}
\newtheorem{corollary}[corollary]{Corollary}
\aliascntresetthe{corollary}

\theoremstyle{definition}
\newaliascnt{definition}{theorem}
\newtheorem{definition}[definition]{Definition}
\aliascntresetthe{definition}
\newaliascnt{assumption}{theorem}
\newtheorem{assumption}[assumption]{Assumption}
\aliascntresetthe{assumption}
\newaliascnt{example}{theorem}
\newtheorem{example}[example]{Example}
\aliascntresetthe{example}

\theoremstyle{remark}
\newaliascnt{remark}{theorem}
\newtheorem{remark}[remark]{Remark}
\aliascntresetthe{remark}

\newcommand{\R}{\mathbb{R}}
\newcommand{\N}{\mathbb{N}}
\newcommand{\Tend}{\mathfrak T}
\newcommand{\HH}{\mathcal H}
\newcommand{\VV}{\mathcal V}
\newcommand{\W}{\mathscr W}
\newcommand{\norm}[1]{\left\lVert #1 \right\rVert}
\newcommand{\inner}[2]{\left( #1,#2 \right)}
\newcommand{\dual}[2]{\left\langle #1,#2 \right\rangle}
\newcommand{\diff}{\mathop{}\!\mathrm{d}}

\title{Information Limits for Stable Reduction of Coupled Diffusion Systems}
\author{Hiroki Ishizaka\\
Team FEM, Matsuyama, Japan\\
E-mail: \texttt{h.ishizaka005@gmail.com}}
\date{}

\begin{document}
\maketitle

\begin{abstract}
Eliminating an unresolved diffusion field from a coupled parabolic system produces an operator-valued memory term and a contribution from the hidden initial state. We study whether a reduced prediction of the retained field meets a prescribed tolerance. For known hidden initial data, the discarded response, damped by diffusion and transmitted through the coupling, gives an upper error bound; an auxiliary diffusion problem supplies a lower bound. Rigorous discretisation bounds then yield tests that accept, reject or remain undecided. If the discarded initial state is known only through a norm bound, the minimum worst-case error equals that bound multiplied by the time-uniform gain from the missing state to the retained trajectory. Our main result encloses this infinite-dimensional gain using a finite spectral Galerkin evolution, its residual in the full equations and an explicit diffusion tail. The enclosure converges without compactness or smallness of the coupling. Strict comparisons of the optimal radius with a tolerance have finite certificates if all endpoint quantities admit arbitrarily accurate rigorous enclosures. Multiplication-coupling examples with infinitely many dynamically coupled modes illustrate both settings, including explicit finite-element error budgets.
\end{abstract}



\section{Introduction}
Large systems of evolution equations often contain more variables than can be resolved in a practical computation. A reduced model retains selected variables and approximates the influence of the rest. Discarded variables may continue to affect the retained dynamics, so an efficient reduction need not give an accurate prediction. In a linear system, eliminating unresolved variables produces a memory term and a forcing contribution from their initial states and sources. For coupled diffusion, the memory is operator-valued and can be approximated by retaining finitely many hidden modes. We ask whether the resulting prediction meets a prescribed tolerance, and what accuracy can be guaranteed when part of the hidden initial state is unavailable.

We consider the passive system
\begin{align*}
  \partial_t u+Au+C^*v&=f,\\
  \partial_t v+Lv-Cu&=g,
\end{align*}
on $H=L^2(\Omega)$ for a bounded domain $\Omega$, with the precise assumptions given in Section~\ref{sec:full-system}. The diffusion operators $A$ and $L$ are coercive and self-adjoint with compact resolvents, and $C$ is bounded but need not be compact. Eliminating $v$ gives the kernel $K(t)=C^*e^{-tL}C$. We use the spectral damping of $L$ to control the influence of discarded modes on the retained trajectory. When the full initial state is prescribed, the question concerns the error of a specified reduced model. When part of the hidden initial state is missing, the available information may instead impose a positive lower limit on the accuracy that any prediction can guarantee uniformly over the admissible states.

Fix a finite horizon $\Tend>0$. Let $P_N$ project onto the first $N$ hidden diffusion modes, ordered by nondecreasing eigenvalue, and set $Q_N=I-P_N$. Suppose that $h=Q_Nv_0$ is known only through $\norm{h}_H\le R$, where $v_0$ is the hidden initial state and $R\ge0$; all other data are fixed. Linearity and the symmetry of this ball imply that the minimum worst-case error in $\mathcal C([0,\Tend];H)$ is $R\Gamma_{N,\Tend}$ (Theorem~\ref{thm:information-barrier}). Here, $\Gamma_{N,\Tend}$ is the operator norm of the map from $h$ to its contribution to the retained trajectory. This is the radius of information in optimal recovery \cite{MicchelliRivlin1977,TraubWasilkowskiWozniakowski1988}. The minimax identity is elementary; the difficulty is to enclose the gain of the infinite-dimensional evolution. Our main result, Theorem~\ref{thm:finite-gain-enclosure}, gives
\begin{align*}
  (\gamma-\Delta)_+
  \le \Gamma_{N,\Tend}
  \le
  \bigl((\gamma+\Delta)^2
  +\overline\Gamma_{N',\Tend}^{\,2}\bigr)^{1/2},
  \quad N<N',
\end{align*}
where $(a)_+=\max\{a,0\}$. The quantity $\gamma$ is the gain of a finite spectral Galerkin system on the initial modes $N+1,\ldots,N'$, $\Delta$ bounds its residual in the full equations, and $\overline\Gamma_{N',\Tend}$ bounds the contribution of modes beyond $N'$ by diffusion damping. To obtain a prescribed enclosure width, first choose $N'$ sufficiently large and then enlarge the Galerkin spaces for this fixed $N'$. No compactness or smallness of $C$ is required. If the spectral, coupling and finite-evolution quantities entering the endpoints admit arbitrarily accurate rigorous enclosures, every strict comparison of $R\Gamma_{N,\Tend}$ with a tolerance has a finite certificate (Corollary~\ref{cor:finite-termination}). This is a conditional existence result, and no decision is asserted at equality. The bounds use finite evolutions and estimates of their full-space residuals, without computing the full coupled trajectory.

For fully known hidden initial data, Theorem~\ref{thm:resolved-certificate} bounds the error of the $N$-mode reduction by
\begin{align*}
  \norm{u-u_N}_{\mathcal C([0,\Tend];H)}
  \le
  \norm{C^*a_N}_{L^1(0,\Tend;H)},
\end{align*}
where $a_N$ is the discarded response generated by $Q_Nv_0$ and $Q_N(Cu_N+g)$. The resulting explicit estimate is at most $\eta_N/\lambda_{N+1}$ times the standard product-energy residual bound, where $\eta_N=\norm{Q_NC}_{\mathcal L(H)}$ and $\lambda_{N+1}$ is the first discarded eigenvalue. This prefactor is less than one when $\eta_N<\lambda_{N+1}$. An auxiliary diffusion problem gives a lower bound (Theorem~\ref{thm:two-sided-retained-decision}), and Theorem~\ref{thm:discrete-three-way} includes numerical error. The test accepts, rejects or remains undecided according to the resulting error interval. An undecided outcome means that the bounds do not settle the prescribed tolerance; it is not a diagnosis of instability. Corollary~\ref{cor:robust-missing-initial-data} combines these known-data bounds with the gain enclosure to certify a tolerance for all admissible missing states or to show that some admissible state violates it. The estimates require rigorous bounds for the relevant spectral data, coupling terms and numerical evolutions.

Memory terms produced by elimination go back to Zwanzig and Mori \cite{Zwanzig1961,Mori1965}; see \cite{ZhuDominyVenturi2018} for estimates of the memory integral. Conservative realisations of dissipative memory are studied in \cite{FigotinSchenker2007}, and passivity-preserving reduction in \cite{BreitenUnger2022}. Residual-based a posteriori bounds for reduced parabolic problems go back at least to Grepl and Patera \cite{GreplPatera2005}. Rettberg et al.\ \cite{RettbergEtAl2024} develop hierarchical bounds and bounds based on an auxiliary linear problem for port-Hamiltonian systems. Our auxiliary diffusion response follows this error-estimation viewpoint. Model reduction with nonzero initial conditions is studied in \cite{HeinkenschlossReisAntoulas2011,BeattieGugercinMehrmann2017,SchroederVoigt2023,LiljegrenSailer2024}. In particular, Liljegren-Sailer \cite{LiljegrenSailer2024} uses observability Gramians and an offline--online splitting to estimate the initial-condition contribution to the $L^2$-output error of stable finite-dimensional systems. We instead enclose a time-uniform gain on an infinite-dimensional initial-state space, expressed by a supremum of pointwise observation-operator norms rather than a single time-integrated Gramian (Remark~\ref{rem:pointwise-gramian}). Error bounds for infinite-dimensional systems are classical \cite{GloverCurtainPartington1988,CurtainZwart2020,BeckerHartmann2019}. Lazar and Zuazua \cite{LazarZuazua2024} and Grubi\v{s}i\'c and Lazar \cite{GrubisicLazar2026} derive eigenvalue bounds for Gramians, the latter allowing non-compact control operators under its stated hypotheses. These spectral estimates concern a different quantity from the finite-horizon time-uniform gain enclosed here.

Sections~\ref{sec:full-system}--\ref{sec:stable-reduction} introduce the coupled system, its exact reduction, passive storage and kernel-perturbation stability. Much of this framework appears in the author's earlier note \cite{IshizakaStableReduction2026}; we give the diffusion arguments to keep the paper self-contained. Proposition~5.13 of that note establishes finite-mode passive truncation for a bounded hidden generator under a compact-coupling assumption. Here, $\lambda_{N+1}\to\infty$ yields convergence in the time-integrated operator norm for every bounded $C$ (Theorem~\ref{thm:spectral-kernel}). Section~\ref{sec:spectral-reduction} develops spectral reduction and the error bound for known data. Section~\ref{sec:information-limit} establishes the information limit and its convergent enclosure. Section~\ref{sec:two-sided} develops two-sided decisions and their discretisation. Multiplication-coupling examples illustrate both settings without restricting the full dynamics to finitely many modes. Non-invariant truncations and nonlinear hidden dynamics remain outside the present analysis.

\section{Coupled resolved--unresolved diffusion system}
\label{sec:full-system}
Let $\Omega\subset\R^d$, $d\ge1$, be a bounded Lipschitz domain.  We  define real Hilbert spaces as
\begin{align*}
  H&:=L^2(\Omega),
  \quad V:=H^1_0(\Omega) \quad \text{with } \norm{w}_V:=\norm{\nabla w}_{L^2(\Omega)^d}.
\end{align*}
From the Poincar\'e inequality, this is a norm on $V$, and
\begin{align*}
  V\hookrightarrow H\cong H^*\hookrightarrow V^*
\end{align*}
is a Gelfand triple with continuous and dense embeddings.  For a finite time horizon $\Tend>0$, we use the standard evolution space
\begin{align}\label{eq:evolution-space}
  \W(0,\Tend)
  :=\bigl\{w\in L^2(0,\Tend;V):
  \partial_t w\in L^2(0,\Tend;V^*)\bigr\}.
\end{align}
It is continuously embedded into $\mathcal C([0,\Tend];H)$; see, for example, Lions and Magenes \cite[Chapter~III]{LionsMagenes1972}.

We consider two bounded symmetric bilinear forms
\begin{align*}
  a_A,a_L:V\times V\to\R.
\end{align*}
The first describes the diffusion of the resolved field and the second that of the unresolved field.

\begin{assumption}\label{ass:forms}
There are constants $M_A,M_L>0$ and $\alpha_A,\alpha_L>0$ such that
\begin{align}
  |a_A(w,\varphi)|&\le M_A\norm{w}_V\norm{\varphi}_V,
  \quad a_A(w,w)\ge \alpha_A\norm{w}_V^2,
  \label{eq:form-A}\\
  |a_L(w,\varphi)|&\le M_L\norm{w}_V\norm{\varphi}_V,
  \quad a_L(w,w) \ge \alpha_L\norm{w}_V^2
  \label{eq:form-L}
\end{align}
for any $w,\varphi\in V$.
\end{assumption}

\begin{remark}[Role of the standing assumptions]\label{rem:assumption-role}
Not every result uses all of these hypotheses.  Well-posedness and the energy estimate in \cref{thm:full-wp} need only boundedness and coercivity of the block form together with the skew coupling.  Symmetry of $a_L$ makes the memory kernel self-adjoint and positive, and it is needed for the spectral calculus later on.  On a bounded domain with Dirichlet conditions, the self-adjoint realisations automatically have compact resolvent, since $V\hookrightarrow H$ is compact; this is not used in \cref{sec:elimination,sec:passive-memory,sec:stable-reduction} and enters only from \cref{sec:spectral-reduction} on, through eigenmode projections. Coercivity of $a_L$ and Poincar\'e's inequality give the spectral gap used in the time-integrated kernel estimates, so no separate gap condition is imposed.
\end{remark}

Let $A,L:V\to V^*$ be the operators induced by the forms,
\begin{align*}
  \dual{Aw}{\varphi}:=a_A(w,\varphi),
  \quad
  \dual{Lw}{\varphi}:=a_L(w,\varphi).
\end{align*}
Their positive self-adjoint realisations in $H$ are characterised by
\begin{align*}
  D(A)&:=\bigl\{w\in V:\ \exists h\in H\ \text{ such that }a_A(w,\varphi)=\inner{h}{\varphi}_H
  \ \forall\varphi\in V\bigr\},\\
  D(L)&:=\bigl\{w\in V:\ \exists h\in H\ \text{ such that }a_L(w,\varphi)=\inner{h}{\varphi}_H
  \ \forall\varphi\in V\bigr\},
\end{align*}
with $Aw=h$ and $Lw=h$, respectively.  As usual, we use the same symbols $A$ and $L$ for the form operators $V\to V^*$ and for these self-adjoint realisations.  The meaning is determined by the expression: in variational equations $Aw$ and $Lw$ are elements of $V^*$, whereas $e^{-tA}$, $e^{-tL}$, $D(L^\theta)$, the spectral measure $E_L$, and all spectral projectors refer to the realisation in $H$.  On $D(A)$ and $D(L)$ the two meanings agree through the embedding $H\hookrightarrow V^*$.  We write $\dual{\cdot}{\cdot}$ for the $V^*,V$-duality pairing and $\inner{\cdot}{\cdot}_H$ for the pivot inner product; under the canonical identification $H\subset V^*$, these pairings coincide whenever the first argument belongs to $H$.

Let
\begin{align*}
  C\in\mathcal L(H)
\end{align*}
be a bounded coupling operator, with Hilbert-space adjoint $C^*$.

The resolved variable is denoted by $u$, while $v$ denotes the unresolved variable. We assume
\begin{align}\label{eq:data-full}
  u_0,v_0\in H,
  \quad
  f,g\in L^2(0,\Tend;H).
\end{align}
The variational theory would allow right-hand sides in $L^2(0,\Tend;V^*)$.  We take them in $L^2(0,\Tend;H)$ because the hidden variation-of-constants formula, the operator-valued convolution and the estimates of the hidden contribution are all formulated in $H$.

We consider
\begin{subequations}\label{eq:full-system}
\begin{align}
  \partial_tu(t)+Au(t)+C^*v(t)&=f(t),
  \label{eq:full-u}\\
  \partial_tv(t)+Lv(t)-Cu(t)&=g(t),
  \label{eq:full-v}
\end{align}
\end{subequations}
in $V^*$, with
\begin{align}\label{eq:full-initial}
  u(0)=u_0,
  \quad
  v(0)=v_0.
\end{align}
With these signs, the off-diagonal operator is skew-adjoint on $H\times H$, so its contribution to the energy balance vanishes.

\subsection{Variational formulation and well-posedness}
We set
\begin{align*}
  \HH:=H\times H,
  \quad
  \VV:=V\times V,
\end{align*}
with the natural product norms.  For $U=(u,v)$ and $\Phi=(\varphi,\psi)$ in $\VV$, we define
\begin{align}\label{eq:block-form}
  \mathfrak a(U,\Phi)
  := {}&a_A(u,\varphi)+a_L(v,\psi)
  +\inner{C^*v}{\varphi}_H
  -\inner{Cu}{\psi}_H.
\end{align}
Because the spaces are real,
\begin{align}\label{eq:coupling-cancel}
  \inner{C^*v}{u}_H-\inner{Cu}{v}_H=0.
\end{align}
Consequently,
\begin{align}\label{eq:block-coercivity}
  \mathfrak a(U,U)
  =a_A(u,u)+a_L(v,v)
  \ge \alpha\norm{U}_{\VV}^2,
  \quad
  \alpha:=\min\{\alpha_A,\alpha_L\}.
\end{align}
Notice that the coercivity constant in \eqref{eq:block-coercivity} does not contain $\norm{C}_{\mathcal L(H)}$.

\begin{theorem}[Well-posedness of the coupled system]\label{thm:full-wp}
Let \cref{ass:forms} hold and let the data satisfy \eqref{eq:data-full}. Then, problem \eqref{eq:full-system}--\eqref{eq:full-initial} has a unique solution
\begin{align*}
  (u,v)\in \W(0,\Tend)\times\W(0,\Tend).
\end{align*}
Furthermore, for $t\in[0,\Tend]$,
\begin{align}\label{eq:full-energy}
  &\frac12\norm{u(t)}_H^2
  +\frac12\norm{v(t)}_H^2
  +\int_0^t a_A(u(s),u(s))\,\diff s
  +\int_0^t a_L(v(s),v(s))\,\diff s
  \notag\\
  &\quad
  =\frac12\norm{u_0}_H^2
  +\frac12\norm{v_0}_H^2
  +\int_0^t\inner{f(s)}{u(s)}_H\,\diff s
  +\int_0^t\inner{g(s)}{v(s)}_H\,\diff s.
\end{align}
In particular, with $F:=(f,g)$ and $U_0:=(u_0,v_0)$,
\begin{align}\label{eq:full-stability}
  \sup_{0\le t\le\Tend}\norm{(u(t),v(t))}_{\HH}^2
  +\alpha\norm{(u,v)}_{L^2(0,\Tend;\VV)}^2
  \le
  2\norm{U_0}_{\HH}^2
  +2\alpha^{-1}\norm{F}_{L^2(0,\Tend;\VV^*)}^2.
\end{align}
The constant on the right-hand side is independent of the size of the coupling operator $C$.
\end{theorem}

\begin{proof}
The form $\mathfrak a$ is bounded on $\VV\times\VV$.  Indeed, the only point not contained in \cref{ass:forms} is the coupling term, for which the continuous embedding $V\hookrightarrow H$ gives
\begin{align*}
  |\inner{C^*v}{\varphi}_H|
  +|\inner{Cu}{\psi}_H|
  \le c\norm{C}_{\mathcal L(H)}
  \norm{(u,v)}_{\VV}\norm{(\varphi,\psi)}_{\VV}.
\end{align*}
The coercivity is exactly \eqref{eq:block-coercivity}.  The standard theorem for coercive variational evolution equations therefore gives a unique solution in $\W(0,\Tend)\times\W(0,\Tend)$; see \cite[Chapter~III]{LionsMagenes1972}.

For the energy identity, take $\varphi=u(t)$ in \eqref{eq:full-u} and $\psi=v(t)$ in \eqref{eq:full-v}, add the two equalities, and use \eqref{eq:coupling-cancel}.  The chain rule in the Gelfand triple yields
\begin{align*}
  \dual{\partial_tu}{u}
  =\frac12\frac{\diff}{\diff t}\norm{u}_H^2,
  \quad
  \dual{\partial_tv}{v}
  =\frac12\frac{\diff}{\diff t}\norm{v}_H^2
\end{align*}
for almost every $t$.  Integration in time gives \eqref{eq:full-energy}.

Finally, use \eqref{eq:block-coercivity} and Young's inequality in
\begin{align*}
  \dual{F}{U}_{\VV^*,\VV}
  \le \norm{F}_{\VV^*}\norm{U}_{\VV}
  \le \frac{1}{2\alpha}\norm{F}_{\VV^*}^2
  +\frac{\alpha}{2}\norm{U}_{\VV}^2.
\end{align*}
It follows, for every $t\in[0,\Tend]$, that
\begin{align*}
  \norm{(u(t),v(t))}_{\HH}^2
  +\alpha\norm{(u,v)}_{L^2(0,t;\VV)}^2
  \le
  \norm{U_0}_{\HH}^2
  +\alpha^{-1}\norm{F}_{L^2(0,\Tend;\VV^*)}^2.
\end{align*}
Taking the supremum over $t$ in the first term and $t=\Tend$ in the second, and adding, gives \eqref{eq:full-stability}.
\end{proof}

\begin{remark}[The skew coupling]\label{rem:passive-coupling}
The cancellation \eqref{eq:coupling-cancel} is what allows the hidden field to be eliminated without losing the dissipative structure of the full system.  It reappears in \cref{sec:passive-memory} as a storage identity for the memory term.
\end{remark}

\subsection{A concrete diffusion model}
A typical example is the following. Let
\begin{align}\label{eq:concrete-forms}
\begin{aligned}
  a_A(w,\varphi)
  &:=\nu\int_\Omega \nabla w\cdot\nabla\varphi\,\diff x
  +a\int_\Omega w\varphi\,\diff x,\\
  a_L(w,\varphi)
  &:=\kappa\int_\Omega \nabla w\cdot\nabla\varphi\,\diff x
  +\ell\int_\Omega w\varphi\,\diff x,
\end{aligned}
\end{align}
where
\begin{align*}
  \nu>0,
  \quad
  \kappa>0,
  \quad
  a,\ell\ge0.
\end{align*}
Let $c\in L^\infty(\Omega;\R)$ and we define
\begin{align*}
  (Cw)(x):=c(x)w(x).
\end{align*}
Then, $C=C^*$ on $H$, and \eqref{eq:full-system} becomes
\begin{subequations}\label{eq:concrete-pde}
\begin{align}
  \partial_tu-\nu\Delta u+au+c(x)v&=f,
  \label{eq:concrete-u}\\
  \partial_tv-\kappa\Delta v+\ell v-c(x)u&=g
  \label{eq:concrete-v}
\end{align}
\end{subequations}
in $\Omega\times(0,\Tend)$, with homogeneous Dirichlet conditions for both fields.  A constant coupling $c(x)\equiv c_0$ gives $C=c_0I$.

\section{Exact elimination of the unresolved field}
\label{sec:elimination}
The $H$-realisation of $L$ is positive and self-adjoint. Let
\begin{align*}
  S_L(t):=e^{-tL},
  \quad t\ge0,
\end{align*}
be the contraction semigroup generated by $-L$.  We set
\begin{align}\label{eq:L-gap}
  \lambda_L:=\inf \sigma(L)>0.
\end{align}
where $\sigma(L)$ denotes the spectrum of the self-adjoint realisation of $L$ in $H$. The strict positivity follows from the coercivity of $a_L$ and Poincar\'e's inequality.  Equivalently, $L\ge\lambda_L I$ on $H$.  Therefore,
\begin{align}\label{eq:SL-decay}
  \norm{S_L(t)}_{\mathcal L(H)}\le e^{-\lambda_L t}
  \quad(t\ge0).
\end{align}
The semigroup representation used below is standard; see, for example, Pazy \cite[Chapter~4]{Pazy1983} or Engel and Nagel \cite[Chapter~II]{EngelNagel2000}.

\begin{lemma}[Variation of constants for the hidden field]\label{lem:voc-v}
Let $(u,v)$ be the solution of \cref{thm:full-wp}. Then, for any $t\in[0,\Tend]$,
\begin{align}\label{eq:v-representation}
  v(t)
  =S_L(t)v_0
  +\int_0^t S_L(t-s)\bigl(Cu(s)+g(s)\bigr)\,\diff s
\end{align}
in $H$.
\end{lemma}

\begin{proof}
We set $h:=Cu+g$. Because $u\in L^2(0,\Tend;V)\hookrightarrow L^2(0,\Tend;H)$ and $C\in\mathcal L(H)$, one has $h\in L^2(0,\Tend;H)$.  The positive self-adjoint operator in $H$ associated with $a_L$ generates the contraction semigroup $S_L$.  For $v_0\in H$ and $h\in L^2(0,\Tend;H)$, the function
\begin{align*}
  w(t):=S_L(t)v_0+\int_0^tS_L(t-s)h(s)\,\diff s
\end{align*}
is the mild solution of $w'+Lw=h$.  By the standard correspondence between the semigroup solution and the variational solution associated with the closed coercive form $a_L$, one has $w\in\W(0,\Tend)$ and $w$ is the unique variational solution with initial value $v_0$; see, for example, \cite[Chapter~4]{Pazy1983}. Because $v$ is that variational solution, $v=w$ in $\mathcal C([0,\Tend];H)$, which proves \eqref{eq:v-representation} for any $t\in[0,\Tend]$.
\end{proof}

\begin{definition}[Memory kernel generated by the hidden diffusion]\label{def:kernel}
For $t\ge0$, we define
\begin{align}\label{eq:kernel-def}
  K(t):=C^*S_L(t)C=C^*e^{-tL}C
  \in\mathcal L(H).
\end{align}
For $u\in L^1(0,\Tend;H)$, we set
\begin{align}\label{eq:kernel-convolution}
  (K*u)(t):=\int_0^t K(t-s)u(s)\,\diff s.
\end{align}
\end{definition}

The contribution of the hidden initial state and the hidden forcing is written as
\begin{align}\label{eq:hidden-contribution}
  r_{v_0,g}(t)
  :=C^*S_L(t)v_0
  +\int_0^t C^*S_L(t-s)g(s)\,\diff s.
\end{align}
We then set
\begin{align}\label{eq:reduced-force}
  F_{v_0,g}(t):=f(t)-r_{v_0,g}(t).
\end{align}
Substitution of \eqref{eq:v-representation} into \eqref{eq:full-u} gives the reduced equation
\begin{align}\label{eq:reduced-problem}
  \partial_tu(t)+Au(t)+(K*u)(t)=F_{v_0,g}(t)
  \quad\text{in }V^*,
  \quad
  u(0)=u_0.
\end{align}
The memory term is not a modelling choice here. It is generated by the semigroup of the hidden diffusion.

\begin{proposition}[Elementary properties of the generated kernel]\label{prop:kernel-properties}
The kernel $K$ defined by \eqref{eq:kernel-def} has the following properties.
\begin{enumerate}[label=\textup{(\roman*)},leftmargin=2.4em]
  \item For any $w\in H$, the map $t\mapsto K(t)w$ is continuous on $[0,\infty)$.  As an $\mathcal L(H)$-valued map, $K$ is continuous in operator norm on $(0,\infty)$ and is Bochner measurable on any finite interval.
  \item $K(t)$ is self-adjoint and non-negative for any $t\ge0$.
  \item For any $t\ge0$,
  \begin{align}\label{eq:kernel-decay}
    \norm{K(t)}_{\mathcal L(H)}
    \le \norm{C}_{\mathcal L(H)}^2e^{-\lambda_Lt}.
  \end{align}
  In particular,
  \begin{align}\label{eq:kernel-L1}
    \norm{K}_{L^1(0,\Tend;\mathcal L(H))}
    \le
    \frac{\norm{C}_{\mathcal L(H)}^2}{\lambda_L}
    \bigl(1-e^{-\lambda_L\Tend}\bigr).
  \end{align}
\end{enumerate}
\end{proposition}

\begin{proof}
Strong-operator continuity follows from that of $S_L$.  To see the operator-norm continuity away from the origin, fix $t_0>0$.  For $s,t\ge t_0/2$, the spectral theorem and the mean-value theorem give
\begin{align*}
  \norm{e^{-tL}-e^{-sL}}_{\mathcal L(H)}
  \le |t-s|\sup_{\lambda\ge\lambda_L}\lambda e^{-(t_0/2)\lambda},
\end{align*}
and the supremum on the right is finite. Therefore, $S_L$, and with it $K$, is norm-continuous on $(0,\Tend]$.  A norm-continuous map on $(0,\Tend]$ has separable range, and adding the single value $K(0)$ does not change this, so $K:[0,\Tend]\to\mathcal L(H)$ is Bochner measurable.  Continuity at $t=0$ plays no role here.

Because $L$ is self-adjoint and non-negative, so is $S_L(t)$.  Therefore,
\begin{align*}
  K(t)^*=C^*S_L(t)^*C=K(t)
\end{align*}
and, for $w\in H$,
\begin{align*}
  \inner{K(t)w}{w}_H
  =\inner{S_L(t)Cw}{Cw}_H\ge0.
\end{align*}
Finally, \eqref{eq:SL-decay} gives \eqref{eq:kernel-decay}, and integration in time gives \eqref{eq:kernel-L1}.
\end{proof}

\begin{lemma}[Regularity of the reduced terms]\label{lem:reduced-terms}
Let $K$ be given by \eqref{eq:kernel-def}. Then, for any $u\in L^2(0,\Tend;H)$,
\begin{align}\label{eq:convolution-L2-bound}
  \norm{K*u}_{L^2(0,\Tend;H)}
  \le
  \norm{K}_{L^1(0,\Tend;\mathcal L(H))}
  \norm{u}_{L^2(0,\Tend;H)}.
\end{align}
Furthermore, the hidden contribution \eqref{eq:hidden-contribution} belongs to $L^2(0,\Tend;H)$ and satisfies
\begin{align}\label{eq:hidden-contribution-L2}
  \norm{r_{v_0,g}}_{L^2(0,\Tend;H)}
  \le{}&
  \norm{C}_{\mathcal L(H)}
  \left(\frac{1-e^{-2\lambda_L\Tend}}{2\lambda_L}\right)^{1/2}
  \norm{v_0}_H \notag\\
  &+\frac{\norm{C}_{\mathcal L(H)}}{\lambda_L}
  \bigl(1-e^{-\lambda_L\Tend}\bigr)
  \norm{g}_{L^2(0,\Tend;H)}.
\end{align}
Consequently, $F_{v_0,g}\in L^2(0,\Tend;H)$, and the reduced equation \eqref{eq:reduced-problem} is well defined in $L^2(0,\Tend;V^*)$ for any $u\in\W(0,\Tend)$.
\end{lemma}

\begin{proof}
Estimate \eqref{eq:convolution-L2-bound} is Young's inequality for Banach-valued convolution.  For the first term in \eqref{eq:hidden-contribution}, \eqref{eq:SL-decay} gives
\begin{align*}
  \norm{C^*S_L(\cdot)v_0}_{L^2(0,\Tend;H)}
  \le
  \norm{C}_{\mathcal L(H)}
  \left(\int_0^{\Tend}e^{-2\lambda_Lt}\,\diff t\right)^{1/2}
  \norm{v_0}_H.
\end{align*}
For the second term, Young's inequality and \eqref{eq:SL-decay} yield
\begin{align*}
  \norm{\int_0^{\cdot}C^*S_L(\cdot-s)g(s)\,\diff s}_{L^2(0,\Tend;H)}
  \le
  \norm{C}_{\mathcal L(H)}
  \norm{e^{-\lambda_L\cdot}}_{L^1(0,\Tend)}
  \norm{g}_{L^2(0,\Tend;H)},
\end{align*}
which is \eqref{eq:hidden-contribution-L2}.  Finally, $\W(0,\Tend)\hookrightarrow L^2(0,\Tend;H)$, so the convolution term is $H$-valued in $L^2$, while $Au\in L^2(0,\Tend;V^*)$.
\end{proof}

\begin{theorem}[Exact reduction and reconstruction]\label{thm:equivalence}
Let \cref{ass:forms} hold and let the data satisfy \eqref{eq:data-full}.
\begin{enumerate}[label=\textup{(\roman*)},leftmargin=2.4em]
  \item If $(u,v)$ solves \eqref{eq:full-system}--\eqref{eq:full-initial}, then $u$ solves the reduced problem \eqref{eq:reduced-problem}.
  \item Conversely, let $u\in\W(0,\Tend)$ solve \eqref{eq:reduced-problem}, and we define $v$ as \eqref{eq:v-representation}. Then, $v\in\W(0,\Tend)$ and $(u,v)$ solves \eqref{eq:full-system}--\eqref{eq:full-initial}.
\end{enumerate}
Consequently, the full and reduced formulations are equivalent, and the unresolved field is uniquely reconstructed from the resolved trajectory and the hidden data.
\end{theorem}

\begin{proof}
By \cref{lem:reduced-terms}, all terms in the reduced equation belong to the claimed spaces.  For part (i), \cref{lem:voc-v} gives \eqref{eq:v-representation}.  Applying $C^*$ and using the boundedness of $C^*$ to pass it through the Bochner integrals gives, for any $t$,
\begin{align*}
  C^*v(t)=r_{v_0,g}(t)+(K*u)(t).
\end{align*}
Substitution into \eqref{eq:full-u} proves \eqref{eq:reduced-problem}.

For part (ii), we put $h:=Cu+g\in L^2(0,\Tend;H)$ and we define $v$ by \eqref{eq:v-representation}.  The argument in \cref{lem:voc-v} shows that $v\in\W(0,\Tend)$ and that
\begin{align*}
  \partial_tv+Lv=Cu+g
  \quad\text{in }L^2(0,\Tend;V^*).
\end{align*}
Furthermore, strong continuity of $S_L$ and the vanishing of the convolution term as $t\downarrow0$ give $v(0)=v_0$ in $H$.  Therefore, \eqref{eq:full-v}--\eqref{eq:full-initial} hold.  The same identity $C^*v=r_{v_0,g}+K*u$, now read in $L^2(0,\Tend;H)$, converts the reduced equation exactly into \eqref{eq:full-u}. Therefore, $(u,v)$ solves the full system.  If two reduced solutions existed, reconstructing their hidden fields would give two solutions of the full problem with the same data, contradicting \cref{thm:full-wp}; hence both reconstruction and the reduced trajectory are unique.
\end{proof}

\begin{corollary}[Well-posedness of the generated memory equation]\label{cor:reduced-wp}
Under the assumptions of \cref{thm:equivalence}, problem \eqref{eq:reduced-problem} has a unique solution in $\W(0,\Tend)$.
\end{corollary}

\begin{remark}\label{rem:volterra-context}
Equation \eqref{eq:reduced-problem} is a linear Volterra evolution equation with an operator-valued kernel, and the general theory of such equations is classical; see Pr\"uss \cite{Pruss1993}. What we use below is the particular structure inherited from the hidden diffusion: positivity of the kernel and the associated storage identity.
\end{remark}

\section{Passive memory and the storage identity}
\label{sec:passive-memory}

\begin{definition}[Positive-type operator kernel]\label{def:positive-type}
For $G\in L^1(0,\Tend;\mathcal L(H))$ and $w\in L^2(0,\Tend;H)$, the causal convolution $(G*w)(t):=\int_0^tG(t-s)w(s)\,\diff s$ is defined for almost every $t$ and belongs to $L^2(0,\Tend;H)$ by Young's inequality. We say that $G$ is of positive type on $[0,\Tend]$ if
\begin{align}\label{eq:positive-type-def}
  \int_0^{\Tend}\inner{(G*w)(t)}{w(t)}_H\,\diff t\ge0
  \quad
  \forall w\in L^2(0,\Tend;H).
\end{align}
\end{definition}

\begin{lemma}[Localisation of positive type]\label{lem:positive-type-local}
If $G$ is of positive type on $[0,\Tend]$, then for any $t\in[0,\Tend]$,
\begin{align}\label{eq:positive-type-local}
  \int_0^t\inner{(G*w)(s)}{w(s)}_H\,\diff s\ge0
  \quad\forall w\in L^2(0,t;H).
\end{align}
\end{lemma}

\begin{proof}
We extend $w$ by zero to $(0,\Tend)$. Its convolution is unchanged before $t$, and the extended test function vanishes afterwards. Applying \eqref{eq:positive-type-def} gives \eqref{eq:positive-type-local}; the case $t=0$ is immediate.
\end{proof}

For $G\in L^1(0,\Tend;\mathcal L(H))$, $w\in L^2(0,\Tend;H)$ and $t\in[0,\Tend]$, we write
\begin{align}\label{eq:memory-work-functional}
  \mathscr D_G[w;t]
  :=\int_0^t\inner{(G*w)(s)}{w(s)}_H\,\diff s
\end{align}
for the accumulated work of the memory term. We call the memory passive when this work is non-negative for every input $w$ and terminal time $t$, with zero initial memory state. By \cref{lem:positive-type-local}, positive-type kernels have this property.

For $u\in L^2(0,\Tend;H)$, we define the hidden response with zero initial state and zero hidden forcing as
\begin{align}\label{eq:zu-def}
  z_u(t):=\int_0^t S_L(t-s)Cu(s)\,\diff s.
\end{align}
As in \cref{lem:voc-v}, $z_u\in\W(0,\Tend)$ is the unique variational
solution of
\begin{align}\label{eq:zu-eq}
  \partial_tz_u+Lz_u=Cu,
  \quad
  z_u(0)=0,
\end{align}
with the evolution equation in $L^2(0,\Tend;V^*)$ and the initial condition in $H$. Passing the bounded operator $C^*$ through the Bochner integral gives
\begin{align}\label{eq:Ku-Cz}
  (K*u)(t)=C^*z_u(t)
\end{align}
for almost every $t\in(0,\Tend)$, with both sides in $L^2(0,\Tend;H)$.

\begin{theorem}[Storage--dissipation identity]\label{thm:storage}
For any $u\in L^2(0,\Tend;H)$ and any $t\in[0,\Tend]$,
\begin{align}\label{eq:storage-identity}
  \int_0^t\inner{(K*u)(s)}{u(s)}_H\,\diff s
  ={}&\frac12\norm{z_u(t)}_H^2
  +\int_0^t a_L(z_u(s),z_u(s))\,\diff s
  \ge0.
\end{align}
In particular, $K$ is of positive type on any finite time interval.
\end{theorem}

\begin{proof}
By \eqref{eq:Ku-Cz} and \eqref{eq:zu-eq},
\begin{align*}
  \inner{(K*u)(s)}{u(s)}_H
  &=\inner{z_u(s)}{Cu(s)}_H\\
  &=\dual{\partial_tz_u(s)}{z_u(s)}
    +a_L(z_u(s),z_u(s))
\end{align*}
for almost every $s$. Integrating, using the Gelfand-triple chain rule and $z_u(0)=0$, gives \eqref{eq:storage-identity}. The continuous $H$-valued representative of $z_u$ makes the identity valid for any $t\in[0,\Tend]$.
\end{proof}

In \eqref{eq:storage-identity}, the terminal term is the stored hidden energy and the integral of $a_L$ is the accumulated diffusion dissipation.

\begin{corollary}[Energy identity for the reduced system]\label{cor:reduced-energy}
Under the standing assumptions, suppose $v_0=0$ and $g=0$. Let $u\in\W(0,\Tend)$ solve
\begin{align}\label{eq:reduced-zero-hidden}
  \partial_tu+Au+K*u=f,
  \quad u(0)=u_0,
\end{align}
and let $z_u$ be defined by \eqref{eq:zu-def}. Then, for any $t\in[0,\Tend]$,
\begin{align}\label{eq:reduced-energy}
  &\frac12\norm{u(t)}_H^2
  +\frac12\norm{z_u(t)}_H^2
  +\int_0^t a_A(u(s),u(s))\,\diff s
  +\int_0^t a_L(z_u(s),z_u(s))\,\diff s
  \notag\\
  &\quad
  =\frac12\norm{u_0}_H^2
  +\int_0^t\inner{f(s)}{u(s)}_H\,\diff s.
\end{align}
\end{corollary}

\begin{proof}
Test \eqref{eq:reduced-zero-hidden} by $u$, integrate, and use \cref{thm:storage}.
\end{proof}

\textit{Spectral form of the memory kernel.} Let $E_L$ be the projection-valued spectral measure of the positive self-adjoint $H$-realisation of $L$. Thus, $E_L(B)$ is an orthogonal spectral projection for each Borel set $B\subset[\lambda_L,\infty)$, and $E_L([\lambda_L,\infty))$ is the identity on $H$. All functions and powers of $L$ below are defined by the spectral calculus. We set
\begin{align}\label{eq:Sigma-def}
  \Sigma(B):=C^*E_L(B)C.
\end{align}
Each $\Sigma(B)$ is self-adjoint and non-negative. For $x,y\in H$, the scalar set function $\mu_{x,y}(B):=\inner{\Sigma(B)x}{y}_H =\inner{E_L(B)Cx}{Cy}_H$ is a finite signed Borel measure. This is the meaning of weak countable additivity of $\Sigma$; also, $\Sigma([\lambda_L,\infty))=C^*C$. The spectral theorem gives
\begin{align}\label{eq:kernel-spectral}
  K(t)
  =\int_{[\lambda_L,\infty)}e^{-\lambda t}\,\diff\Sigma(\lambda),
  \quad t\ge0,
\end{align}
in the weak operator sense, namely $\inner{K(t)x}{y}_H =\int_{[\lambda_L,\infty)}e^{-\lambda t}\,\diff\mu_{x,y}(\lambda)$ for all $x,y\in H$.

For bounded self-adjoint operators, the Loewner order $B_1\preceq B_2$ means $\inner{(B_2-B_1)w}{w}_H\ge0$ for any $w\in H$. A function $G:(0,\infty)\to\mathcal L(H)$ with self-adjoint values is completely monotone in this order if it is infinitely differentiable in operator norm and $0\preceq(-1)^mG^{(m)}(t)$ for $m=0,1,\ldots$ and $t>0$, where $G^{(m)}$ is the $m$th time derivative and $G^{(0)}=G$. For scalar functions, the same derivative inequalities use the usual order on $\mathbb R$.

\begin{proposition}[Operator-valued completely monotone structure]\label{prop:cm-structure}
The map $K:(0,\infty)\to\mathcal L(H)$ is infinitely differentiable in operator norm, and
\begin{align}\label{eq:operator-cm-derivatives}
  (-1)^m K^{(m)}(t)
  =C^*L^m e^{-tL}C
  \quad(m=0,1,\ldots,\ t>0).
\end{align}
The operator on the right is self-adjoint and non-negative, so $K$ is completely monotone in the Loewner order. Equivalently, for any $w\in H$, $k_w(t):=\inner{K(t)w}{w}_H$ is completely monotone. Indeed, the non-negative Borel measure $\mu_w(B):=\norm{E_L(B)Cw}_H^2$ has total mass $\norm{Cw}_H^2$, and
\begin{align}
  k_w(t)
 &=\int_{[\lambda_L,\infty)}e^{-\lambda t}\,\diff\mu_w(\lambda), \label{eq:scalar-spectral}\\
 (-1)^m k_w^{(m)}(t)
  &=\int_{[\lambda_L,\infty)}\lambda^m e^{-\lambda t}\,\diff\mu_w(\lambda)
  \ge0 \label{eq:cm-derivatives}
\end{align}
for any $m=0,1,\ldots$ and any $t>0$.
\end{proposition}

\begin{proof}
For $q\ge0$ and $t>0$, the spectral theorem gives $e^{-tL}H\subset D(L^q)$ and
\begin{align*}
  \norm{L^q e^{-tL}}_{\mathcal L(H)}
  \le\sup_{\lambda\ge\lambda_L}\lambda^q e^{-t\lambda}<\infty.
\end{align*}
On each compact positive-time interval, the multipliers $\lambda^j e^{-t\lambda}$ are uniformly bounded for any $j=0,1,\ldots$. Taylor's formula for the scalar exponential therefore justifies repeated differentiation in operator norm, giving
\begin{align*}
  \frac{\diff^m}{\diff t^m}e^{-tL}=(-1)^mL^m e^{-tL}.
\end{align*}
This proves \eqref{eq:operator-cm-derivatives}. Furthermore,
\begin{align*}
  C^*L^m e^{-tL}C
  =\bigl(L^{m/2}e^{-tL/2}C\bigr)^*
   \bigl(L^{m/2}e^{-tL/2}C\bigr),
\end{align*}
with bounded factors, proving self-adjointness and non-negativity. Formula \eqref{eq:scalar-spectral} follows from \eqref{eq:kernel-spectral}, and pairing \eqref{eq:operator-cm-derivatives} with $w$ gives \eqref{eq:cm-derivatives}.
\end{proof}

\begin{remark}[Direction of the construction]\label{rem:direction-memory}
The construction in \cite{IshizakaCoercivity2026,IshizakaGraphSpace2026} starts with a scalar completely monotone kernel. Its Bernstein measure, the non-negative Borel measure whose Laplace transform is that kernel, supplies internal variables for the energy and well-posedness analysis. Here, the hidden dynamics produces $\Sigma$ through \eqref{eq:Sigma-def} and hence $K(t)=C^*e^{-tL}C$. In general, this kernel does not factor into a scalar time kernel and a fixed spatial operator, so spectral truncation can change both its temporal dependence and its spatial action.
\end{remark}

\section{Stable reduction from unresolved influence to retained trajectories}
\label{sec:stable-reduction}
We compare retained trajectories when the kernel and the hidden-data contribution are approximated.

Let $S_A(t):=e^{-tA}$ be the semigroup of the self-adjoint $H$-realisation of $A$, with $\norm{S_A(t)}_{\mathcal L(H)}\le1$. For $G\in L^1(0,\Tend;\mathcal L(H))$, $w_0\in H$ and $q\in L^1(0,\Tend;H)$, we consider
\begin{align}\label{eq:general-reduced}
  \partial_tw+Aw+G*w=q,\quad w(0)=w_0.
\end{align}
A mild solution is a function $w\in\mathcal C([0,\Tend];H)$ satisfying
\begin{align}\label{eq:general-mild}
  w(t)=S_A(t)w_0
  +\int_0^t S_A(t-s)\bigl(q(s)-(G*w)(s)\bigr)\,\diff s
\end{align}
for any $t\in[0,\Tend]$. This Volterra equation has a unique mild solution; see \cite{Pruss1993}. For trajectories in $\W(0,\Tend)$, $Aw$ denotes the form operator in $V^*$, as in \cref{sec:full-system}. The energy estimates below assume that these trajectories exist and require only $a_A(w,w)\ge0$.

\begin{definition}[Stable reduction in the $L^1$ topology]\label{def:L1-consistent-reduction}
Let $(K,r)$ and $(K_N,r_N)$, $N\in\N$, be the kernel and hidden-data contribution of the full and reduced models, respectively, with $K,K_N\in L^1(0,\Tend;\mathcal L(H))$ and $r,r_N\in L^1(0,\Tend;H)$. We call the family an \emph{$L^1$-consistent stable reduction} on $[0,\Tend]$ if
\begin{align}\label{eq:stable-reduction-def}
  \norm{K_N-K}_{L^1(0,\Tend;\mathcal L(H))}
  +\norm{r_N-r}_{L^1(0,\Tend;H)}
  \longrightarrow0
  \quad(N\to\infty).
\end{align}
For a common resolved forcing $f$, the corresponding right-hand sides are $F=f-r$ and $F_N=f-r_N$. When $r_N=r=0$, we call $\{K_N\}$ an \emph{$L^1$-consistent kernel reduction}. The reduction is called \emph{passive} if every $K_N$ is of positive type in the sense of \cref{def:positive-type}. Consistency also makes $\{K_N\}$ bounded in $L^1(0,\Tend;\mathcal L(H))$.
\end{definition}

\begin{proposition}[Volterra kernel-to-trajectory stability]\label{prop:volterra-kernel-trajectory}
Let $u_N,u\in\mathcal C([0,\Tend];H)$ be the mild solutions of
\begin{align}
  \partial_tu_N+Au_N+K_N*u_N&=F_N,
  \quad u_N(0)=u_{0,N},
  \label{eq:reduced-N}\\
  \partial_tu+Au+K*u&=F,
  \quad u(0)=u_0,
  \label{eq:reduced-limit}
\end{align}
where $u_{0,N},u_0\in H$. The kernels satisfy $K_N,K\in L^1(0,\Tend;\mathcal L(H))$, and the forcing terms satisfy $F_N,F\in L^1(0,\Tend;H)$. We set $M_N:=\norm{K_N}_{L^1(0,\Tend;\mathcal L(H))}$. Then,
\begin{align}\label{eq:volterra-kernel-trajectory-estimate}
  \norm{u_N-u}_{\mathcal C([0,\Tend];H)}
  \le{}& e^{M_N\Tend}
  \Bigl(
  \norm{u_{0,N}-u_0}_H
  +\norm{F_N-F}_{L^1(0,\Tend;H)}
  \notag\\
  &\hspace{27mm}
  +\norm{K_N-K}_{L^1(0,\Tend;\mathcal L(H))}
   \norm{u}_{L^1(0,\Tend;H)}
  \Bigr).
\end{align}
Thus, for $L^1$-consistent kernels, convergent data give $u_N\to u$ in $\mathcal C([0,\Tend];H)$.
\end{proposition}

\begin{proof}
We set $e_N:=u_N-u$ and $R_N:=F_N-F-(K_N-K)*u$. Subtracting the mild equations gives
\begin{align*}
  e_N(t)=S_A(t)e_N(0)
  +\int_0^t S_A(t-s)\bigl(R_N(s)-(K_N*e_N)(s)\bigr)\,\diff s.
\end{align*}
Contractivity and Fubini's theorem yield
\begin{align*}
  \norm{e_N(t)}_H
  \le \norm{e_N(0)}_H+\norm{R_N}_{L^1(0,\Tend;H)}
       +M_N\int_0^t\norm{e_N(s)}_H\,\diff s.
\end{align*}
Young's convolution inequality bounds the residual by
\begin{align*}
  \norm{R_N}_{L^1(0,\Tend;H)}
  \le \norm{F_N-F}_{L^1(0,\Tend;H)}
      +\norm{K_N-K}_{L^1(0,\Tend;\mathcal L(H))}
       \norm{u}_{L^1(0,\Tend;H)}.
\end{align*}
Gronwall's inequality proves \eqref{eq:volterra-kernel-trajectory-estimate}. The convergence follows because consistency bounds $M_N$ uniformly.
\end{proof}

\begin{theorem}[Passivity-enhanced kernel-to-trajectory stability]\label{thm:kernel-trajectory}
Let $A:V\to V^*$ be induced by a bounded bilinear form $a_A:V\times V\to\R$ satisfying
\begin{align}\label{eq:passive-resolved-form}
  a_A(w,w)\ge0\quad\forall w\in V.
\end{align}
Let $u_N,u\in\W(0,\Tend)$ solve \eqref{eq:reduced-N} and \eqref{eq:reduced-limit}, with $K_N,K\in L^1(0,\Tend;\mathcal L(H))$ and $F_N,F\in L^2(0,\Tend;H)$. Suppose $K_N$ is of positive type on $[0,\Tend]$, and we define
\begin{align}\label{eq:passive-residual-size}
  B_N:={}&\norm{F_N-F}_{L^1(0,\Tend;H)}
  +\norm{K_N-K}_{L^1(0,\Tend;\mathcal L(H))}
   \norm{u}_{L^1(0,\Tend;H)}.
\end{align}
Then,
\begin{align}\label{eq:kernel-trajectory-estimate}
  \norm{u_N-u}_{\mathcal C([0,\Tend];H)}
  &\le B_N+\sqrt{B_N^2+\norm{u_{0,N}-u_0}_H^2}
  \notag\\
  &\le \norm{u_{0,N}-u_0}_H+2B_N.
\end{align}
In particular, if $\{K_N\}$ is a passive $L^1$-consistent kernel reduction, $u_{0,N}\to u_0$ in $H$ and $F_N\to F$ in $L^1(0,\Tend;H)$, then
\begin{align}\label{eq:trajectory-convergence}
  u_N\longrightarrow u
  \quad\text{in }\mathcal C([0,\Tend];H).
\end{align}
The estimate has no exponential amplification factor.
\end{theorem}

\begin{proof}
We set $e_N:=u_N-u$ and $R_N:=F_N-F-(K_N-K)*u$. Then,
\begin{align}\label{eq:passive-error-equation}
  \partial_te_N+Ae_N+K_N*e_N=R_N,
  \quad e_N(0)=u_{0,N}-u_0.
\end{align}
Young's inequality gives $R_N\in L^2(0,\Tend;H)$ and $\norm{R_N}_{L^1(0,\Tend;H)}\le B_N$. Testing by $e_N$ and integrating, the non-negative resolved and memory contributions give
\begin{align}\label{eq:passive-error-energy}
  \frac12\norm{e_N(t)}_H^2
  \le\frac12\norm{e_N(0)}_H^2
     +\int_0^t\inner{R_N(s)}{e_N(s)}_H\,\diff s
\end{align}
for any $t\in[0,\Tend]$; see \cref{lem:positive-type-local}. With $X_N:=\norm{e_N}_{\mathcal C([0,\Tend];H)}$, this yields $X_N^2\le\norm{e_N(0)}_H^2+2B_NX_N$. The non-negative root gives the first bound in \eqref{eq:kernel-trajectory-estimate}; the second uses $\sqrt{a^2+b^2}\le a+b$ for $a,b\ge0$.
\end{proof}

\begin{corollary}[Unit gain for a realised passive kernel]\label{cor:realised-passive-unit-gain}
Assume the trajectory, data and resolved-form hypotheses of \cref{thm:kernel-trajectory}. In place of abstract positive type, suppose a finite-dimensional real Hilbert space $\mathcal Y_N$, a non-negative self-adjoint $\widehat L_N\in\mathcal L(\mathcal Y_N)$ and $\widehat C_N\in\mathcal L(H,\mathcal Y_N)$ give
\begin{align}\label{eq:finite-passive-realisation}
  K_N(t)=\widehat C_N^*e^{-t\widehat L_N}\widehat C_N,
  \quad t\ge0.
\end{align}
Here, $\mathcal L(H,\mathcal Y_N)$ denotes bounded linear maps from $H$ to $\mathcal Y_N$, and $\widehat C_N^*$ is the Hilbert-space adjoint. With $B_N$ from \eqref{eq:passive-residual-size},
\begin{align}\label{eq:realised-passive-unit-gain}
  \norm{u_N-u}_{\mathcal C([0,\Tend];H)}
  \le\norm{u_{0,N}-u_0}_H+B_N.
\end{align}
When $u_{0,N}=u_0$, this becomes
\begin{align}\label{eq:realised-passive-unit-gain-same-initial}
  \norm{u_N-u}_{\mathcal C([0,\Tend];H)}\le B_N.
\end{align}
\end{corollary}

\begin{proof}
We set $e_N:=u_N-u$, $R_N:=F_N-F-(K_N-K)*u$ and
\begin{align}\label{eq:realised-passive-error-state}
  \zeta_N(t)
  :=\int_0^t e^{-(t-s)\widehat L_N}\widehat C_Ne_N(s)\,\diff s.
\end{align}
This auxiliary state is absolutely continuous with square-integrable derivative and satisfies $\partial_t\zeta_N+\widehat L_N\zeta_N=\widehat C_Ne_N$, $\zeta_N(0)=0$, and $K_N*e_N=\widehat C_N^*\zeta_N$. Testing the error equation by $e_N$ and the auxiliary equation by $\zeta_N$, the coupling terms cancel:
\begin{align}\label{eq:realised-passive-pointwise-energy}
  &\frac12\frac{\diff}{\diff t}
   \left(\norm{e_N(t)}_H^2+\norm{\zeta_N(t)}_{\mathcal Y_N}^2\right)
  +a_A(e_N(t),e_N(t))
  +\inner{\widehat L_N\zeta_N(t)}{\zeta_N(t)}_{\mathcal Y_N}
  \notag\\
  &\quad=\inner{R_N(t)}{e_N(t)}_H
\end{align}
for almost every $t$. Both dissipation terms are non-negative. We set $E_N(t):=(\norm{e_N(t)}_H^2+ \norm{\zeta_N(t)}_{\mathcal Y_N}^2)^{1/2}$. For $\varepsilon>0$, the chain rule gives
\begin{align*}
  \frac{\diff}{\diff t}(E_N^2+\varepsilon^2)^{1/2}
  \le \norm{R_N(t)}_H
     \frac{\norm{e_N(t)}_H}{(E_N^2+\varepsilon^2)^{1/2}}
  \le\norm{R_N(t)}_H
\end{align*}
almost everywhere. Integration and $\varepsilon\downarrow0$ yield $E_N(t)\le\norm{e_N(0)}_H+\norm{R_N}_{L^1(0,\Tend;H)}$. Because $\norm{e_N(t)}_H\le E_N(t)$ and $\norm{R_N}_{L^1(0,\Tend;H)}\le B_N$, this proves \eqref{eq:realised-passive-unit-gain}. The coefficient one is sharp in this non-negative-form class: take $A=0$, $K_N=K=0$, $e_N(0)=0$ and $R_N(t)=h$ for fixed $0\ne h\in V$. Then, $e_N(t)=th$, and equality holds at $t=\Tend$.
\end{proof}

\begin{remark}[Sign assumptions]\label{rem:passive-minimality}
For the given $\W$-trajectories, the energy estimate uses only $a_A(w,w)\ge0$ and positive type of $K_N$. Without these assumptions, the bound can fail. For $H=V=\R$, initial mismatch one and zero forcing, $A=0$ and $K_N=K\equiv-a$, $a>0$, give $e(t)=\cosh(\sqrt a\,t)$; with $K_N=K=0$ and $a_A(x,y)=-axy$, one gets $e(t)=e^{at}$. We make no claim that the factor two in \eqref{eq:kernel-trajectory-estimate} is optimal. The explicit realisation in \cref{cor:realised-passive-unit-gain} gives coefficient one.
\end{remark}

\begin{remark}[Choice of kernel topology]\label{rem:L1-topology}
The comparison uses the integrated kernel norm in \eqref{eq:stable-reduction-def}. For spectral truncations, this error can tend to zero even when the kernel error at $t=0$ does not; see \cref{sec:spectral-reduction}.
\end{remark}


\section{Spectral reduction of the hidden diffusion}
\label{sec:spectral-reduction}
We retain the forms and data of \cref{sec:full-system}. The operator $L$ has compact resolvent, meaning that $(L+I)^{-1}:H\to H$ is compact: coercivity makes this resolvent bounded from $H$ to $V$, and $V\hookrightarrow H$ is compact. Here, $I$ denotes the identity on $H$. No eigenfunction expansion was needed in the preceding sections. We choose an orthonormal eigenbasis $(\phi_j)_{j\ge1}$ of $H$, with eigenvalues listed with multiplicity:
\begin{align}\label{eq:eigenpairs}
  L\phi_j=\lambda_j\phi_j,\quad
  0<\lambda_1\le\lambda_2\le\cdots,\quad\lambda_j\to\infty.
\end{align}
Thus, $\lambda_1=\lambda_L$ from \eqref{eq:L-gap}. Let $H_N:=\operatorname{span}\{\phi_1,\ldots,\phi_N\}\subset D(L)$, equipped with the inner product inherited from $H$, let $P_N$ be the orthogonal projection in $H$ onto $H_N$, and set $Q_N:=I-P_N$. Then, $LP_N=P_NL$ on $D(L)$, and $P_N$ commutes with $e^{-tL}$. We define
\begin{align}\label{eq:LN-CN}
  L_N:=L|_{H_N}:H_N\to H_N,\quad C_N:=P_NC:H\to H_N.
\end{align}
The operator $L_N$ is positive and self-adjoint, and $C_N^*z=C^*z$ for $z\in H_N\subset H$. The generated kernel is
\begin{align}\label{eq:KN-def}
  K_N(t):=C_N^*e^{-tL_N}C_N=C^*P_Ne^{-tL}P_NC,\quad t\ge0.
\end{align}
It is self-adjoint, non-negative and of positive type by the storage identity. The coupling to the discarded modes is measured by
\begin{align}\label{eq:etaN-def}
  \eta_N:=\norm{Q_NC}_{\mathcal L(H)}\le\norm{C}_{\mathcal L(H)}.
\end{align}

\subsection{Exact finite reduction and generated kernel structure}
Only the hidden field is truncated. We consider
\begin{subequations}\label{eq:finite-full-system}
\begin{align}
  \partial_tu_N+Au_N+C_N^*v_N&=f,
  \label{eq:finite-full-u}\\
  \partial_tv_N+L_Nv_N-C_Nu_N&=P_Ng,
  \label{eq:finite-full-v}
\end{align}
\end{subequations}
with
\begin{align}\label{eq:finite-full-initial}
  u_N(0)=u_0,\quad v_N(0)=P_Nv_0.
\end{align}
The equations hold almost everywhere in $V^*$ and $H_N$, respectively, using the Gelfand triple $V\times H_N\hookrightarrow H\times H_N\hookrightarrow V^*\times H_N$. The coupling cancels in the product energy, since $\inner{C_N^*z}{w}_H=\inner{z}{C_Nw}_{H_N}$.

\begin{theorem}[Exact finite reduction and commutation]\label{thm:finite-reduction-commute}
For any $N\in\N$, \eqref{eq:finite-full-system}--\eqref{eq:finite-full-initial} has a unique solution $u_N\in\W(0,\Tend)$, $v_N\in H^1(0,\Tend;H_N)$. Here, $H^1(0,\Tend;H_N)$ consists of absolutely continuous paths whose values and derivatives belong to $L^2(0,\Tend;H_N)$. For any $t\in[0,\Tend]$,
\begin{align}\label{eq:finite-full-energy}
  &\frac12\norm{u_N(t)}_H^2+\frac12\norm{v_N(t)}_{H_N}^2
  +\int_0^t\bigl[a_A(u_N,u_N)
                  +\inner{L_Nv_N}{v_N}_{H_N}\bigr]\,\diff s
  \notag\\
  &\quad=\frac12\norm{u_0}_H^2+\frac12\norm{P_Nv_0}_{H_N}^2
  +\int_0^t\bigl[\inner{f}{u_N}_H
                 +\inner{P_Ng}{v_N}_{H_N}\bigr]\,\diff s.
\end{align}
For $w\in L^2(0,\Tend;H)$, we define
\begin{align}
  \mathscr H[w](t)&:=e^{-tL}v_0
    +\int_0^t e^{-(t-s)L}(Cw(s)+g(s))\,\diff s,
    \label{eq:full-hidden-response-map}\\
  \mathscr H_N[w](t)&:=e^{-tL_N}P_Nv_0
    +\int_0^t e^{-(t-s)L_N}(C_Nw(s)+P_Ng(s))\,\diff s.
    \label{eq:finite-hidden-response-map}
\end{align}
Then, for any $t\in[0,\Tend]$,
\begin{align}\label{eq:hidden-response-commutation}
  \mathscr H_N[w](t)=P_N\mathscr H[w](t).
\end{align}
In particular,
\begin{align}\label{eq:vN-representation}
  v_N(t)=\mathscr H_N[u_N](t),
\end{align}
and elimination gives
\begin{align}\label{eq:general-spectral-reduced}
  \partial_tu_N+Au_N+K_N*u_N=F_N,\quad u_N(0)=u_0,
\end{align}
where
\begin{align}\label{eq:FN-def}
  F_N&:=f-r_N,\quad
  r_N(t):=C^*P_Ne^{-tL}v_0
    +\int_0^t C^*P_Ne^{-(t-s)L}g(s)\,\diff s.
\end{align}
Equivalently,
\begin{align}\label{eq:rN-finite-realisation}
  r_N(t)=C_N^*e^{-tL_N}P_Nv_0
    +\int_0^t C_N^*e^{-(t-s)L_N}P_Ng(s)\,\diff s.
\end{align}
Conversely, any variational solution of \eqref{eq:general-spectral-reduced}, reconstructed by \eqref{eq:vN-representation}, solves the truncated coupled system. Furthermore,
\begin{align}
  C^*\mathscr H[w]&=r_{v_0,g}+K*w,
    \label{eq:full-response-elimination}\\
  C_N^*\mathscr H_N[w]&=r_N+K_N*w.
    \label{eq:finite-response-elimination}
\end{align}
Thus, truncation and elimination commute, including the hidden data.
\end{theorem}

\begin{proof}
The block form on $V\times H_N$ is bounded, and its diagonal is at least $\alpha_A\norm{w}_V^2+\lambda_1\norm{z}_{H_N}^2$. Variational theory gives the unique solution. The hidden equation then gives $v_N\in H^1(0,\Tend;H_N)$. Testing with $(u_N,v_N)$ gives \eqref{eq:finite-full-energy}.

The identity $e^{-tL_N}P_N=P_Ne^{-tL}$ proves \eqref{eq:hidden-response-commutation} and \eqref{eq:rN-finite-realisation}. Applying $C^*$ and $C_N^*$ to the response formulas gives \eqref{eq:full-response-elimination}--\eqref{eq:finite-response-elimination}. With $w=u_N$, substitution yields the reduced equation. Conversely, variation of constants gives the hidden equation and initial value, while $C_N^*v_N=r_N+K_N*u_N$ recovers the resolved equation.
\end{proof}

\begin{remark}[The hidden contribution has to be truncated too]\label{rem:kernel-only-not-commute}
For non-zero hidden data, the truncated system uses $(K_N,r_N)$, with $F_N-F_{v_0,g}=r_{v_0,g}-r_N$. Keeping $r_{v_0,g}$ while replacing only $K$ need not describe this system. For $v_0=g=0$, both contributions vanish.
\end{remark}

\begin{theorem}[Coupling-tail spectral estimate and compactness criterion]\label{thm:spectral-kernel}
For any $N\in\N$, $t\ge0$ and $x\in H$,
\begin{align}\label{eq:spectral-kernel-factor}
  K(t)-K_N(t)=(Q_NC)^*e^{-tL}(Q_NC),
\end{align}
\begin{align}\label{eq:spectral-kernel-quadratic}
  \inner{(K(t)-K_N(t))x}{x}_H
    =\norm{e^{-tL/2}Q_NCx}_H^2.
\end{align}
Consequently,
\begin{align}\label{eq:spectral-kernel-pointwise}
  \norm{K(t)-K_N(t)}_{\mathcal L(H)}
    \le\eta_N^2e^{-\lambda_{N+1}t},
\end{align}
\begin{align}\label{eq:spectral-kernel-L1}
  \norm{K-K_N}_{L^1(0,\Tend;\mathcal L(H))}
    \le\frac{\eta_N^2}{\lambda_{N+1}}
             (1-e^{-\lambda_{N+1}\Tend}),
\end{align}
\begin{align}\label{eq:spectral-kernel-uniform}
  \sup_{0\le t\le\Tend}\norm{K(t)-K_N(t)}_{\mathcal L(H)}=\eta_N^2,
\end{align}
and, for $0<\delta\le\Tend$,
\begin{align}\label{eq:spectral-kernel-away-zero}
  \sup_{\delta\le t\le\Tend}\norm{K(t)-K_N(t)}_{\mathcal L(H)}
    \le\eta_N^2e^{-\lambda_{N+1}\delta}.
\end{align}
The kernels are $L^1$-consistent on every finite interval. The following conditions are equivalent:
\begin{enumerate}[label=\textup{(\alph*)},leftmargin=2.4em]
  \item $K_N\to K$ uniformly in operator norm on $[0,\Tend]$
        for one, and hence any, $\Tend>0$;
  \item $\eta_N\to0$;
  \item $C$ is compact on $H$;
  \item $K\in\mathcal C([0,\Tend];\mathcal L(H))$
        for one, and hence any, $\Tend>0$.
\end{enumerate}
\end{theorem}

\begin{proof}
Commutation and $P_NQ_N=0$ give $e^{-tL}-P_Ne^{-tL}P_N=Q_Ne^{-tL}Q_N$, which proves \eqref{eq:spectral-kernel-factor} and \eqref{eq:spectral-kernel-quadratic}. The discarded spectrum gives $\norm{e^{-tL}Q_N}_{\mathcal L(H)}=e^{-\lambda_{N+1}t}$; hence \eqref{eq:spectral-kernel-pointwise}, its time integral \eqref{eq:spectral-kernel-L1}, and \eqref{eq:spectral-kernel-away-zero}. At $t=0$, the difference is $(Q_NC)^*(Q_NC)$, whose norm is $\eta_N^2$, proving \eqref{eq:spectral-kernel-uniform}. Because $\eta_N\le\norm{C}_{\mathcal L(H)}$ and $\lambda_{N+1}\to\infty$, the kernels are $L^1$-consistent.

The uniform identity gives (a)$\Leftrightarrow$(b). If $C$ is compact, the strongly convergent projections $P_N\to I$ converge uniformly on the image under $C$ of the unit ball, so $\eta_N\to0$. Conversely, $\norm{C-P_NC}_{\mathcal L(H)}=\eta_N\to0$ makes $C$ a norm limit of finite-rank operators. Thus, (b)$\Leftrightarrow$(c). Each $K_N$ is continuous in operator norm, so (a) implies (d). For the converse, spectral calculus gives, for $t>0$,
\begin{align*}
  (1-e^{-\lambda_{N+1}t})\norm{Q_NCx}_H^2
    \le\inner{(K(0)-K(t))x}{x}_H.
\end{align*}
Taking the supremum over $\norm{x}_H\le1$, then $N\to\infty$ and $t\downarrow0$, gives $\eta_N\to0$. This proves (d)$\Rightarrow$(b).
\end{proof}

\begin{remark}[Repeated eigenvalues and the choice of truncation]\label{rem:repeated-eigenvalues}
Even if $\lambda_N=\lambda_{N+1}$, the space $H_N$ is invariant because $L$ acts as a scalar on each eigenspace. Thus, no spectral gap at the truncation index is needed. Cutting a repeated eigenspace can make $H_N$ and $K_N$ depend on its chosen basis. For a basis-independent truncation, we use $P_\Lambda:=E_L([\lambda_1,\Lambda])$, $\Lambda\ge\lambda_1$, with $E_L$ from \cref{sec:passive-memory}; the same estimates use the first discarded eigenvalue in place of $\lambda_{N+1}$.
\end{remark}

\begin{theorem}[Monotone passive truncation and discarded-mode balance]\label{thm:discarded-balance}
For bounded self-adjoint operators, write $B_1\preceq B_2$ when $\inner{(B_2-B_1)x}{x}_H\ge0$ for every $x\in H$. Then, for any $t\ge0$,
\begin{align}\label{eq:kernel-loewner}
  0\preceq K_N(t)\preceq K_{N+1}(t)\preceq K(t).
\end{align}
For $w\in L^2(0,\Tend;H)$, let
\begin{align}\label{eq:discarded-state}
  z_{N,w}^{\perp}(t):=\int_0^t e^{-(t-s)L}Q_NCw(s)\,\diff s.
\end{align}
This belongs to $\W(0,\Tend)$, takes values in $Q_NH$, and solves
\begin{align}\label{eq:discarded-state-evolution}
  \partial_tz_{N,w}^{\perp}+Lz_{N,w}^{\perp}=Q_NCw,
  \quad z_{N,w}^{\perp}(0)=0.
\end{align}
For any $t\in[0,\Tend]$, the memory work in \eqref{eq:memory-work-functional} satisfies
\begin{align}\label{eq:discarded-balance}
  \mathscr D_K[w;t]-\mathscr D_{K_N}[w;t]
  =\frac12\norm{z_{N,w}^{\perp}(t)}_H^2
   +\int_0^t a_L(z_{N,w}^{\perp},z_{N,w}^{\perp})\,\diff s\ge0,
\end{align}
\begin{align}\label{eq:discarded-balance-rate}
  0\le\mathscr D_K[w;t]-\mathscr D_{K_N}[w;t]
  \le\frac{\eta_N^2}{\lambda_{N+1}}(1-e^{-\lambda_{N+1}t})
          \norm{w}_{L^2(0,t;H)}^2.
\end{align}
Furthermore,
\begin{align}\label{eq:memory-work-monotone}
  0\le\mathscr D_{K_N}[w;t]\le\mathscr D_{K_{N+1}}[w;t]
    \le\mathscr D_K[w;t],
\end{align}
and $\mathscr D_{K_N}[w;\,\cdot\,]\to\mathscr D_K[w;\,\cdot\,]$ uniformly on $[0,\Tend]$.
\end{theorem}

\begin{proof}
We set $R_N:=P_{N+1}-P_N$. Orthogonality and commutation give
\begin{align*}
  K_{N+1}(t)-K_N(t)=(R_NC)^*e^{-tL}(R_NC).
\end{align*}
Together with \eqref{eq:spectral-kernel-factor}, this proves \eqref{eq:kernel-loewner}. The variational--semigroup correspondence gives \eqref{eq:discarded-state-evolution}, and commutation gives $Q_Nz_{N,w}^{\perp}=z_{N,w}^{\perp}$. Applying \cref{thm:storage} with coupling $Q_NC$ gives \eqref{eq:discarded-balance}; applying it with $R_NC$ gives \eqref{eq:memory-work-monotone}. This uses positive type of the difference kernels, which does not follow from pointwise Loewner order alone. Finally, Young's inequality and Cauchy--Schwarz give
\begin{align*}
  \mathscr D_K[w;t]-\mathscr D_{K_N}[w;t]
  \le\norm{K-K_N}_{L^1(0,t;\mathcal L(H))}
       \norm{w}_{L^2(0,t;H)}^2.
\end{align*}
The kernel bound proves \eqref{eq:discarded-balance-rate}. It also gives uniform convergence: the error is bounded by $\lambda_{N+1}^{-1}\norm{C}_{\mathcal L(H)}^2$ times $\norm{w}_{L^2(0,\Tend;H)}^2$.
\end{proof}

\begin{remark}[The role of compactness]\label{rem:compactness-role-section6}
Compactness of $C$ characterises uniform kernel convergence up to $t=0$. Positive-time, $L^1$-kernel and uniform memory-work convergence hold for every bounded $C$. The equivalence with continuity of $K$ at zero uses $\lambda_j\to\infty$; it is not asserted for general positive self-adjoint generators without compact resolvent.
\end{remark}

\begin{remark}[Spectral powers]\label{rem:spectral-powers}
For $\theta\in\R$, spectral powers mean
\begin{align*}
  D(L^\theta)&:=\left\{x\in H:
       \sum_{j\ge1}\lambda_j^{2\theta}|\inner{x}{\phi_j}_H|^2<\infty\right\},
  \quad L^\theta x :=\sum_{j\ge1}\lambda_j^\theta\inner{x}{\phi_j}_H\phi_j.
\end{align*}
Negative powers are bounded on $H$ because $\lambda_1>0$. The projections $P_N,Q_N$ commute with these powers on their domains. For $p>0$, the notation $O(\lambda_{N+1}^{-p})$, respectively $o(\lambda_{N+1}^{-p})$, means that $\lambda_{N+1}^p$ times the stated non-negative quantity is bounded, respectively tends to zero, as $N\to\infty$, with the data, operators and $\Tend$ fixed.
\end{remark}

\begin{corollary}[Smoothing coupling and kernel rate]\label{cor:smoothing-coupling-kernel}
If $\theta>0$ and
\begin{align}\label{eq:smoothing-coupling-domain}
  C(H)\subset D(L^\theta),
\end{align}
then $L^\theta C\in\mathcal L(H)$ and
\begin{align}\label{eq:etaN-smoothing}
  \eta_N\le\lambda_{N+1}^{-\theta}\norm{L^\theta C}_{\mathcal L(H)}.
\end{align}
Consequently,
\begin{align}\label{eq:smoothing-kernel-rate}
  \norm{K-K_N}_{L^1(0,\Tend;\mathcal L(H))}
    =O(\lambda_{N+1}^{-1-2\theta}).
\end{align}
The coupling $C$ is compact, and $K_N\to K$ uniformly in operator norm on $[0,\Tend]$.
\end{corollary}

\begin{proof}
Because $C$ is bounded and $L^\theta$ is closed, the everywhere-defined composition $L^\theta C$ has a closed graph, so it is bounded. The factorisation $Q_NC=Q_NL^{-\theta}(L^\theta C)$ and $\norm{Q_NL^{-\theta}}_{\mathcal L(H)}\le\lambda_{N+1}^{-\theta}$ give \eqref{eq:etaN-smoothing}. All conclusions follow from \cref{thm:spectral-kernel}.
\end{proof}

\begin{example}[Compact coupling need not be smoothing]\label{ex:compact-nonsmoothing}
We define $C$ on the eigenbasis as $C\phi_j=c_j\phi_j$, where $c_j:=1/\log(\mathrm e+\lambda_j)$, $\log$ is the natural logarithm and $\mathrm e=\exp(1)$. Then, $C$ is bounded and compact, but $\sup_j\lambda_j^\theta c_j=\infty$ for every $\theta>0$. By \cref{cor:smoothing-coupling-kernel}, $C(H)\not\subset D(L^\theta)$ for any $\theta>0$. Monotonicity of $c_j$ and the diagonal form give
\begin{align*}
  &\eta_N =c_{N+1}, \quad
  \norm{K(t)-K_N(t)}_{\mathcal L(H)}=c_{N+1}^2e^{-\lambda_{N+1}t},\\
  &\norm{K-K_N}_{L^1(0,\Tend;\mathcal L(H))}
    =\frac{1-e^{-\lambda_{N+1}\Tend}}
       {\lambda_{N+1}\log^2(\mathrm e+\lambda_{N+1})}.
\end{align*}
Thus, compactness alone gives no fixed algebraic improvement beyond $\lambda_{N+1}^{-1}$; the smoothing assumption gives such an improvement.
\end{example}

\begin{proposition}[Sharpness of the kernel rate for constant coupling]\label{prop:constant-coupling-sharp}
For $C=c_0I$, $c_0\in\R\setminus\{0\}$, one has $\eta_N=|c_0|$ and
\begin{align}
  \norm{K(t)-K_N(t)}_{\mathcal L(H)}
    &=c_0^2e^{-\lambda_{N+1}t},
    \label{eq:constant-coupling-pointwise-sharp}\\
  \norm{K-K_N}_{L^1(0,\Tend;\mathcal L(H))}
    &=\frac{c_0^2}{\lambda_{N+1}}(1-e^{-\lambda_{N+1}\Tend}).
    \label{eq:constant-coupling-L1-sharp}
\end{align}
In particular,
\begin{align}\label{eq:constant-coupling-asymptotic-sharp}
  \lambda_{N+1}\norm{K-K_N}_{L^1(0,\Tend;\mathcal L(H))}\to c_0^2.
\end{align}
Thus, the kernel rate cannot be improved for general bounded couplings. This does not assert sharpness of the trajectory rate.
\end{proposition}

\begin{proof}
Here, $K(t)-K_N(t)=c_0^2Q_Ne^{-tL}$, whose norm is $c_0^2e^{-\lambda_{N+1}t}$. Integration proves the remaining claims.
\end{proof}

\begin{lemma}[Approximation of the hidden contribution]\label{lem:hidden-contribution-tail}
For $r_{v_0,g}$ from \eqref{eq:hidden-contribution} and $r_N$ from \eqref{eq:FN-def},
\begin{align}\label{eq:hidden-contribution-tail}
  \norm{r_N-r_{v_0,g}}_{L^1(0,\Tend;H)}
  \le\frac{\eta_N}{\lambda_{N+1}}(1-e^{-\lambda_{N+1}\Tend})
       \bigl(\norm{v_0}_H+\norm{g}_{L^1(0,\Tend;H)}\bigr).
\end{align}
This is $O(\lambda_{N+1}^{-1})$, and the bound is $o(\lambda_{N+1}^{-1})$ if $C$ is compact.
\end{lemma}

\begin{proof}
Subtracting the response formulas gives
\begin{align*}
  r_{v_0,g}(t)-r_N(t)=C^*Q_Ne^{-tL}v_0
    +\int_0^t C^*Q_Ne^{-(t-s)L}g(s)\,\diff s.
\end{align*}
Because $\norm{C^*Q_N}_{\mathcal L(H)}=\eta_N$ and $Q_Ne^{-tL}=Q_Ne^{-tL}Q_N$, the operator norm of each propagator is at most $\eta_Ne^{-\lambda_{N+1}t}$. Time integration and Young's $L^1$ inequality prove the estimate. The standing assumption $g\in L^2(0,\Tend;H)$ gives $\norm{g}_{L^1(0,\Tend;H)}\le\Tend^{1/2}\norm{g}_{L^2(0,\Tend;H)}$.
\end{proof}

\begin{theorem}[Resolved trajectory with non-zero hidden data]\label{thm:general-spectral-trajectory}
Let $u$ and $u_N$ be the resolved fields of \eqref{eq:full-system} and \eqref{eq:finite-full-system}, with initial data \eqref{eq:full-initial} and \eqref{eq:finite-full-initial}. Then,
\begin{align}\label{eq:general-spectral-trajectory-rate}
  \norm{u_N-u}_{\mathcal C([0,\Tend];H)}
  \le\frac{1-e^{-\lambda_{N+1}\Tend}}{\lambda_{N+1}}
  \left[\eta_N^2\norm{u}_{L^1(0,\Tend;H)}
  +\eta_N\bigl(\norm{v_0}_H+\norm{g}_{L^1(0,\Tend;H)}\bigr)\right].
\end{align}
Consequently,
\begin{align}\label{eq:general-spectral-trajectory-bigO}
  \norm{u_N-u}_{\mathcal C([0,\Tend];H)}=O(\lambda_{N+1}^{-1}).
\end{align}
No additional spatial regularity of $v_0,g$ is needed. For compact $C$, the displayed bound is $o(\lambda_{N+1}^{-1})$. For fixed data, the trajectory has this improved rate for any bounded $C$; see \cref{cor:fixed-data-little-o}.
\end{theorem}

\begin{proof}
The exact reductions give the kernels $K,K_N$ and forcings $F_{v_0,g},F_N\in L^2(0,\Tend;H)$, with the same resolved initial value. Because $K_N=C_N^*e^{-tL_N}C_N$, apply \cref{cor:realised-passive-unit-gain} and insert \eqref{eq:spectral-kernel-L1} and \eqref{eq:hidden-contribution-tail}.
\end{proof}

\subsection{An error bound for known hidden data}
\label{subsec:known-data-bound}
The preceding estimate involves the exact trajectory $u$. The following bound uses $u_N$ and the prescribed data.

\begin{theorem}[A posteriori bound for the retained field]\label{thm:resolved-certificate}
Under the hypotheses of \cref{thm:general-spectral-trajectory}, let $(u,v)$ and $(u_N,v_N)$ denote the full and truncated solutions. We define
\begin{align}\label{eq:certificate-tail-response}
  a_N(t):=e^{-tL}Q_Nv_0
    +\int_0^t e^{-(t-s)L}Q_N(Cu_N(s)+g(s))\,\diff s.
\end{align}
Then, $a_N\in\W(0,\Tend)$, $a_N(t)\in Q_NH$, and
\begin{align}\label{eq:certificate-basic}
  \norm{u-u_N}_{\mathcal C([0,\Tend];H)}
    \le\norm{C^*a_N}_{L^1(0,\Tend;H)}
    \le\eta_N\norm{a_N}_{L^1(0,\Tend;H)}.
\end{align}
With $\lambda:=\lambda_{N+1}$, we set
\begin{align}\label{eq:certificate-explicit}
  \mathcal E_N:=\frac{\eta_N}{\lambda}\left[
    (1-e^{-\lambda\Tend})\norm{Q_Nv_0}_H
    +\int_0^{\Tend}(1-e^{-\lambda(\Tend-s)})
       \norm{Q_N(Cu_N(s)+g(s))}_H\,\diff s\right].
\end{align}
Then,
\begin{align}\label{eq:certificate-resolved-estimate}
  \norm{u-u_N}_{\mathcal C([0,\Tend];H)}\le\mathcal E_N.
\end{align}
Evaluation requires bounds for the projected data and reduced trajectory, an upper bound for $\eta_N$, and a positive lower bound for $\lambda_{N+1}$. Given these bounds, a rigorously established $\mathcal E_N\le\varepsilon$, for a prescribed tolerance $\varepsilon>0$, guarantees the same trajectory tolerance. Errors in computing $u_N$ are treated in \cref{subsec:fully-discrete-certification}.
\end{theorem}

\begin{proof}
Because $H_N\subset V$, the embedded path $v_N$ belongs to $\W(0,\Tend)$, as does $a_N$ by the variational--semigroup correspondence. We set $e_N:=u-u_N$ and $z_N:=v-v_N-a_N$. Subtracting the equations gives
\begin{align*}
  \partial_te_N+Ae_N+C^*z_N&=-C^*a_N,\quad e_N(0)=0,\\
  \partial_tz_N+Lz_N-Ce_N&=0,\quad z_N(0)=0.
\end{align*}
Testing by $(e_N,z_N)$ and cancelling the coupling gives, almost everywhere,
\begin{align*}
  \frac12\frac{\diff}{\diff t}
    (\norm{e_N}_H^2+\norm{z_N}_H^2)
    +a_A(e_N,e_N)+a_L(z_N,z_N)
    =-\inner{C^*a_N}{e_N}_H.
\end{align*}
Non-negative dissipation and the regularised square-root argument of \cref{cor:realised-passive-unit-gain} yield, for any $t\in[0,\Tend]$,
\begin{align*}
  (\norm{e_N(t)}_H^2+\norm{z_N(t)}_H^2)^{1/2}
    \le\int_0^t\norm{C^*a_N(s)}_H\,\diff s.
\end{align*}
Because $Q_Na_N=a_N$ and $\norm{C^*Q_N}_{\mathcal L(H)}=\eta_N$, this proves \eqref{eq:certificate-basic}. The discarded semigroup estimate gives
\begin{align*}
  \norm{a_N(t)}_H\le e^{-\lambda t}\norm{Q_Nv_0}_H
    +\int_0^t e^{-\lambda(t-s)}
       \norm{Q_N(Cu_N(s)+g(s))}_H\,\diff s.
\end{align*}
Integration and Tonelli's theorem give $\eta_N\norm{a_N}_{L^1(0,\Tend;H)}\le\mathcal E_N$.
\end{proof}

\begin{remark}[Comparison with the standard residual bound]\label{rem:unfiltered-residual-comparison}
The error $(u-u_N,v-v_N)$ has initial value $(0,Q_Nv_0)$ and forcing $(0,Q_N(Cu_N+g))$. The product-energy estimate gives
\begin{align*}
  \norm{u-u_N}_{\mathcal C([0,\Tend];H)}
    \le\norm{Q_Nv_0}_H+\norm{Q_N(Cu_N+g)}_{L^1(0,\Tend;H)}.
\end{align*}
For analogous residual bounds in finite-dimensional passive systems, see \cite[Section~3.1]{RettbergEtAl2024}. In \eqref{eq:certificate-explicit}, the discarded diffusion acts before the coupling, so $\mathcal E_N$ is at most $\eta_N/\lambda_{N+1}$ times this standard bound.
\end{remark}

\begin{corollary}[Resolved-output acceptance criterion]\label{cor:resolved-output-criterion}
Let $\ell\in H^*$, the space of bounded linear functionals on $H$, and let $\varepsilon_\ell>0$. Under the hypotheses of \cref{thm:resolved-certificate},
\begin{align*}
  \sup_{0\le t\le\Tend}|\ell(u(t))-\ell(u_N(t))|
    \le\norm{\ell}_{H^*}\mathcal E_N.
\end{align*}
Thus, $\norm{\ell}_{H^*}\mathcal E_N\le\varepsilon_\ell$ guarantees the output tolerance; otherwise, this test is inconclusive.
\end{corollary}

\begin{proof}
Apply $|\ell(x)|\le\norm{\ell}_{H^*}\norm{x}_H$ and \eqref{eq:certificate-resolved-estimate}.
\end{proof}

\begin{corollary}[Fixed-data rate without compact coupling]\label{cor:fixed-data-little-o}
Under the hypotheses of \cref{thm:general-spectral-trajectory}, for fixed $u_0,v_0,f,g$ and any bounded $C$,
\begin{align*}
  \norm{u-u_N}_{\mathcal C([0,\Tend];H)}=o(\lambda_{N+1}^{-1}).
\end{align*}
This rate need not be uniform on bounded sets of data.
\end{corollary}

\begin{proof}
The baseline estimate gives $u_N\to u$ in $\mathcal C([0,\Tend];H)$. Strong convergence $Q_N\to0$ and dominated convergence give $\norm{Q_Nv_0}_H\to0$, $\norm{Q_Ng}_{L^1(0,\Tend;H)}\to0$, and
\begin{align*}
  \norm{Q_NCu_N}_{L^1(0,\Tend;H)}
    \le\norm{C}_{\mathcal L(H)}\norm{u_N-u}_{L^1(0,\Tend;H)}
       +\norm{Q_NCu}_{L^1(0,\Tend;H)}\to0.
\end{align*}
Because $\eta_N\le\norm{C}_{\mathcal L(H)}$, these limits imply $\lambda_{N+1}\mathcal E_N\to0$ and prove the rate. To see why uniformity can fail, take $A=L$, $C=I$, $u_0=f=g=0$ and $v_0=\phi_{N+1}$. Then, $u_N=0$ and $u(t)=-e^{-\lambda_{N+1}t}\sin(t)\phi_{N+1}$. At $t=\lambda_{N+1}^{-1}\le\Tend$ for large $N$, $\lambda_{N+1}\norm{u(t)}_H\to\mathrm e^{-1}$, although $\norm{v_0}_H=1$.
\end{proof}

\section{Missing hidden initial data}
\label{sec:information-limit}

\subsection{The information limit}
\label{subsec:information-barrier}
Suppose that only $P_Nv_0$ and a bound for $Q_Nv_0$ are available. The missing initial state must then be included in any uniform accuracy guarantee. We identify its optimal worst-case influence, the radius of information in optimal recovery \cite{MicchelliRivlin1977,TraubWasilkowskiWozniakowski1988}.

In this section, write $\mathcal X:=\mathcal C([0,\Tend];H)$ with $\norm{w}_{\mathcal X}:=\max_{0\le t\le\Tend}\norm{w(t)}_H$. Let $\mathscr S(t)$ be the strongly continuous contraction semigroup on $\HH=H\times H$ generated by $\mathcal G U=(-Ax-C^*y,Cx-Ly)^{\top}$, $U=(x,y)^{\top}$, on $D(A)\times D(L)$. Boundedness of the skew coupling and the energy identity in Section~\ref{sec:full-system} give this semigroup; see also \cite{CurtainZwart2020}. With $\pi_u((x,y)^{\top})=x$, we define 
\begin{align}\label{eq:hidden-to-retained-map}
(\mathscr T_{N,\Tend}h)(t):=\pi_u\mathscr S(t)(0,h)^{\top},  \quad h\in Q_NH,\quad 0\le t\le\Tend.
\end{align}
This bounded linear map $Q_NH\to\mathcal X$ has norm $\Gamma_{N,\Tend}:=\norm{\mathscr T_{N,\Tend}}_{\mathcal L(Q_NH,\mathcal X)}$. The subspace $Q_NH$ carries the norm of $H$.

\begin{theorem}[Information limit for discarded initial modes]
\label{thm:information-barrier}
Fix $u_0,f,g,P_Nv_0$ as in Section~\ref{sec:full-system}, $\Tend>0$ and $R\ge0$. Let $u^{(h)}$ be the full retained solution with hidden initial state $P_Nv_0+h$, where $h\in Q_NH$ and $\norm{h}_H\le R$. Then,
\begin{align}\label{eq:unknown-mode-affine}
 u^{(h)}=u^{(0)}+\mathscr T_{N,\Tend}h.
\end{align}
Every prediction $w\in\mathcal X$ independent of $h$ satisfies
\begin{align}\label{eq:information-minimax}
 \sup_{\substack{h\in Q_NH\\\norm{h}_H\le R}}
 \norm{u^{(h)}-w}_{\mathcal X}\ge R\Gamma_{N,\Tend},
\end{align}
with equality for $w=u^{(0)}$. Furthermore,
\begin{align}\label{eq:information-gain-upper}
 0\le\Gamma_{N,\Tend}\le\overline\Gamma_{N,\Tend}
 :=\frac{\eta_N}{\lambda_{N+1}}(1-e^{-\lambda_{N+1}\Tend}),
 \quad \eta_N=\norm{Q_NC}_{\mathcal L(H)},
\end{align}
and
\begin{align}\label{eq:information-zero-characterisation}
 \Gamma_{N,\Tend}=0\quad\Longleftrightarrow\quad C^*Q_N=0
 \quad\Longleftrightarrow\quad Q_NC=0.
\end{align}
Thus, if $R>0$ and $Q_NC\ne0$, no prediction from the known data can guarantee an error strictly below $R\Gamma_{N,\Tend}$ for all admissible $h$.
\end{theorem}

\begin{proof}
Linearity and uniqueness give \eqref{eq:unknown-mode-affine}. Both $h$ and $-h$ are admissible, and
\begin{align*}
 2\norm{\mathscr T_{N,\Tend}h}_{\mathcal X}
 \le\norm{u^{(h)}-w}_{\mathcal X}+\norm{u^{(-h)}-w}_{\mathcal X}.
\end{align*}
Taking the supremum proves \eqref{eq:information-minimax}; setting $w=u^{(0)}$ gives equality, including when $R=0$.

We write $(x,y)^{\top}=\mathscr S(\cdot)(0,h)^{\top}$ and set $a(t)=e^{-tL}h$, $z=y-a$. The zero-initial-state pair $(x,z)^{\top}$ solves $x'+Ax+C^*z=-C^*a$, $z'+Lz-Cx=0$. The energy argument of \cref{thm:resolved-certificate}, using $a(t)\in Q_NH$, gives
\begin{align*}
 \norm{\mathscr T_{N,\Tend}h}_{\mathcal X}
 \le\norm{C^*e^{-tL}h}_{L^1(0,\Tend;H)}
 \le\overline\Gamma_{N,\Tend}\norm{h}_H.
\end{align*}
This proves \eqref{eq:information-gain-upper} and the implication $C^*Q_N=0\Rightarrow\Gamma_{N,\Tend}=0$. Conversely, if $\Gamma_{N,\Tend}=0$, then $x=0$ for any $h\in Q_NH$. The mild first equation is $x(t)=-\int_0^t e^{-(t-s)A}C^*y(s)\,\diff s$. Strong continuity and $y(0)=h$ imply $x(t)/t\to-C^*h$ in $H$ as $t\downarrow0$. Therefore, $C^*Q_N=0$; taking adjoints gives the remaining equivalence.
\end{proof}

The optimal predictor $u^{(0)}$ is a full solution, not a proposed reduced algorithm. The theorem identifies the information limit; when $v_0$ is fully known, \cref{thm:resolved-certificate} instead bounds the error of a specified reduction.

\begin{remark}[Pointwise observation operators and Gramians]
\label{rem:pointwise-gramian}
Let $\mathscr B_N:Q_NH\to\HH$ be the isometry $\mathscr B_Nh=(0,h)^{\top}$ and set
$W_N(t):=(\pi_u\mathscr S(t)\mathscr B_N)^*
          (\pi_u\mathscr S(t)\mathscr B_N)\in\mathcal L(Q_NH)$. These operators are self-adjoint and nonnegative, and
\begin{align}\label{eq:pointwise-gramian-identity}
 \Gamma_{N,\Tend}^2
 =\sup_{0\le t\le\Tend}\norm{W_N(t)}_{\mathcal L(Q_NH)}.
\end{align}
Indeed, the suprema over time and unit initial states commute, and
\[
\norm{\pi_u\mathscr S(t)\mathscr B_N}_{\mathcal L(Q_NH,H)}^2
=\norm{W_N(t)}_{\mathcal L(Q_NH)}.
\]
In contrast, the time-integrated observability Gramian, defined by the bilinear form $\int_0^{\Tend}(W_N(t)h,k)_H\,\diff t$, measures the squared $L^2(0,\Tend;H)$ output norm. Such initial-condition error estimates are studied in \cite{LiljegrenSailer2024}. Identity \eqref{eq:pointwise-gramian-identity} alone gives no finite computation.
\end{remark}

\subsection{Finite-dimensional enclosure of the gain}
\label{subsec:finite-information-gain}
We enclose $\Gamma_{N,\Tend}$ by a finite evolution, its residual in the full equations and a diffusion tail. Numerical use requires rigorous bounds for the spectral data and coupling terms as well as the finite evolution; unverified coupling coefficients do not suffice.

The operator $A$ also has compact resolvent by the assumptions of Section~\ref{sec:full-system}. Fix an orthonormal eigenbasis, ordered by nondecreasing eigenvalue with multiplicity, and let $R_J$ project onto its first $J$ vectors. We put $Q_J^A=I-R_J$ and fix integers $1\le N<N'\le M$, $J\ge1$. For $h\in E_{N,N'}:=(P_{N'}-P_N)H$, solve
\begin{subequations}\label{eq:finite-gain-system}
\begin{align}
 x'_{J,M}+Ax_{J,M}+R_JC^*y_{J,M}&=0,
 \quad x_{J,M}(0) =0,\\
 y'_{J,M}+Ly_{J,M}-P_MCx_{J,M}&=0,
 \quad y_{J,M}(0) =h,
\end{align}
\end{subequations}
in $R_JH\times P_MH$. The two projected coupling terms are adjoints, so this finite system is contractive. Superscript $(k)$ denotes the input $h=\phi_k$. We define
\begin{align}
 \gamma_{J,M;N,N'}
 &:=\sup_{0\le t\le\Tend}
 \norm{h\mapsto x_{J,M}(t;h)}_{\mathcal L(E_{N,N'},H)},
 \label{eq:finite-gain-gamma}\\
 \Delta_{J,M;N,N'}
 &:=\left[\Tend\sum_{k=N+1}^{N'}\int_0^{\Tend}
 \left(\norm{Q_J^AC^*y^{(k)}_{J,M}}_H^2
       +\norm{Q_MCx^{(k)}_{J,M}}_H^2\right)\,\diff t\right]^{1/2}.
 \label{eq:finite-gain-defect}
\end{align}
The sum in \eqref{eq:finite-gain-defect} is independent of the orthonormal basis of $E_{N,N'}$, being the sum of squared norms of the residual map on that basis. The tail bound $\overline\Gamma_{N',\Tend}$ from \eqref{eq:information-gain-upper} tends to zero even with the crude estimate $\eta_{N'}\le\norm{C}_{\mathcal L(H)}$.

\begin{theorem}[Finite-dimensional enclosure of the gain]
\label{thm:finite-gain-enclosure}
For $1\le N<N'\le M$ and $J\ge1$,
\begin{align}\label{eq:finite-gain-enclosure}
 (\gamma_{J,M;N,N'}-\Delta_{J,M;N,N'})_+
 \le\Gamma_{N,\Tend}\le
 \left[(\gamma_{J,M;N,N'}+\Delta_{J,M;N,N'})^2
       +\overline\Gamma_{N',\Tend}^{\,2}\right]^{1/2},
\end{align}
where $(a)_+=\max\{a,0\}$. For each fixed $N'$,
\begin{align}\label{eq:finite-gain-defect-convergence}
 \Delta_{J,M;N,N'}\longrightarrow0
 \quad(J,M\longrightarrow\infty,\ M\ge N').
\end{align}
The enclosure has width at most $2\Delta_{J,M;N,N'}+\overline\Gamma_{N',\Tend}$. For fixed $N$, every desired width is therefore achieved by choosing $N'$ sufficiently large and then $J,M$ sufficiently large for that $N'$. No compactness or smallness of $C$ is required.
\end{theorem}

\begin{proof}
Let $\mathscr T^{J,M}:E_{N,N'}\to\mathcal X$ map $h$ to $x_{J,M}(\cdot;h)$. In the full equations, the finite solution has residual
\begin{align*}
 \mathfrak r_{J,M}(t;h)
 =\bigl(Q_J^AC^*y_{J,M}(t;h),-Q_MCx_{J,M}(t;h)\bigr)^{\top}.
\end{align*}
Its initial state agrees with the full one. The variation-of-constants formula and contractivity give
\begin{align}\label{eq:finite-gain-residual-unit}
 \norm{\mathscr T_{N,\Tend}h-\mathscr T^{J,M}h}_{\mathcal X}
 \le\int_0^{\Tend}\norm{\mathfrak r_{J,M}(t;h)}_{\HH}\,\diff t.
\end{align}
For $h=\sum_{k=N+1}^{N'}h_k\phi_k$, Cauchy--Schwarz gives
\[
\norm{\mathfrak r_{J,M}(t;h)}_{\HH}^2
\le\norm{h}_H^2\sum_{k=N+1}^{N'} \norm{\mathfrak r_{J,M}(t;\phi_k)}_{\HH}^2.
\]
Applying Cauchy--Schwarz also in time yields
\begin{align}\label{eq:finite-gain-low-error}
 \norm{\mathscr T_{N,\Tend}|_{E_{N,N'}}-\mathscr T^{J,M}}
       _{\mathcal L(E_{N,N'},\mathcal X)}\le\Delta_{J,M;N,N'}.
\end{align}
The reverse triangle inequality proves the lower endpoint of \eqref{eq:finite-gain-enclosure}. For arbitrary $h\in Q_NH$, decompose it orthogonally as $h=h_\ell+h_r$, where $h_\ell=(P_{N'}-P_N)h$ and $h_r=Q_{N'}h$. Then,
\begin{align*}
 \norm{\mathscr T_{N,\Tend}h}_{\mathcal X}
 \le(\gamma_{J,M;N,N'}+\Delta_{J,M;N,N'})\norm{h_\ell}_H
     +\overline\Gamma_{N',\Tend}\norm{h_r}_H.
\end{align*}
Cauchy--Schwarz in these two orthogonal components proves the upper endpoint.

For convergence, fix $N'$ and $N<k\le N'$. We write $U_{J,M}^{(k)}=(x_{J,M}^{(k)},y_{J,M}^{(k)})^{\top}$, $U^{(k)}=\mathscr S(\cdot)(0,\phi_k)^{\top}$, $\Pi_{J,M}=R_J\oplus P_M$, $S_0(t)=\operatorname{diag}(e^{-tA},e^{-tL})$ and $\mathbb B((x,y)^{\top})=(C^*y,-Cx)^{\top}$. The initial terms in the two mild equations coincide since $M\ge N'$. As $\Pi_{J,M}$ commutes with $S_0$, subtraction gives
\begin{align*}
 \norm{U_{J,M}^{(k)}(t)-U^{(k)}(t)}_{\HH}
 &\le\norm{C}_{\mathcal L(H)}\int_0^t
 \norm{U_{J,M}^{(k)}(s)-U^{(k)}(s)}_{\HH}\,\diff s
 +\Tend\omega_{J,M}^{(k)},\\
 \omega_{J,M}^{(k)}
 &:=\max_{0\le s\le\Tend}
 \norm{(I-\Pi_{J,M})\mathbb B U^{(k)}(s)}_{\HH}.
\end{align*}
The continuous curve $\mathbb B U^{(k)}([0,\Tend])$ is compact. Strong convergence of the orthogonal projections is uniform on this compact set, so $\omega_{J,M}^{(k)}\to0$. Gronwall's inequality gives $U_{J,M}^{(k)}\to U^{(k)}$ in $\mathcal C([0,\Tend];\HH)$. Furthermore,
\begin{align*}
 \max_{0\le t\le\Tend}\norm{\mathfrak r_{J,M}(t;\phi_k)}_{\HH}
 \le\norm{C}_{\mathcal L(H)}
 \norm{U_{J,M}^{(k)}-U^{(k)}}_{\mathcal C([0,\Tend];\HH)}
 +\omega_{J,M}^{(k)}\longrightarrow0.
\end{align*}
The finite sum in \eqref{eq:finite-gain-defect} now proves \eqref{eq:finite-gain-defect-convergence}. Finally, $\overline\Gamma_{N',\Tend}\le\norm{C}_{\mathcal L(H)}/\lambda_{N'+1}\to0$, and for $\gamma,\Delta,\beta\ge0$, $\sqrt{(\gamma+\Delta)^2+\beta^2}-(\gamma-\Delta)_+\le2\Delta+\beta$. This proves the width statement. No convergence uniform in $N'$, or operator-norm convergence of the full Galerkin semigroups, is asserted.
\end{proof}

\begin{corollary}[Decisions away from the threshold]
\label{cor:finite-termination}
Fix $N$, $\Tend>0$, $R>0$ and $\varepsilon>0$. We denote the endpoints of \eqref{eq:finite-gain-enclosure} by $\gamma^-_{J,M,N'}$ and $\gamma^+_{J,M,N'}$. If $R\Gamma_{N,\Tend}<\varepsilon$, some finite $N',J,M$ satisfy $R\gamma^+_{J,M,N'}<\varepsilon$; if $R\Gamma_{N,\Tend}>\varepsilon$, some satisfy $R\gamma^-_{J,M,N'}>\varepsilon$. Thus, the optimal radius can be compared with the tolerance at a finite truncation level whenever $R\Gamma_{N,\Tend}\ne\varepsilon$, provided the quantities entering the endpoints admit arbitrarily accurate rigorous enclosures. This is an existence statement, not a termination claim for arbitrary operators; nothing is claimed at equality.
\end{corollary}
\begin{proof}
We choose the enclosure width smaller than $|\varepsilon/R-\Gamma_{N,\Tend}|$ using \cref{thm:finite-gain-enclosure}, and enclose the endpoint quantities tightly enough to preserve the resulting strict inequality.
\end{proof}

\begin{corollary}[Uniform tolerance tests with missing initial modes]
\label{cor:robust-missing-initial-data}
For a prediction $w\in\mathcal X$ independent of $h$, we set $E_0=\norm{u^{(0)}-w}_{\mathcal X}$ and $W_N(R;w)=\sup_{h\in Q_NH,\,\norm{h}_H\le R} \norm{u^{(h)}-w}_{\mathcal X}$, where $R\ge0$. If $0\le E_0^-\le E_0\le E_0^+$ and $0\le\gamma^-\le\Gamma_{N,\Tend}\le\gamma^+$, then
\begin{align}\label{eq:robust-missing-tolerance}
 \max\{E_0^-,R\gamma^-\}\le W_N(R;w)\le E_0^++R\gamma^+.
\end{align}
For $\varepsilon>0$, the condition $E_0^++R\gamma^+\le\varepsilon$ guarantees the tolerance for all admissible states, whereas $\max\{E_0^-,R\gamma^-\}>\varepsilon$ implies that at least one state violates it. For $w=u_N$, we use \cref{thm:resolved-certificate} with $Q_Nv_0=0$ to bound $E_0$, and \cref{thm:finite-gain-enclosure} to bound $\Gamma_{N,\Tend}$.
\end{corollary}
\begin{proof}
The triangle inequality in \eqref{eq:unknown-mode-affine} gives the upper bound. The choice $h=0$ and \eqref{eq:information-minimax} give the two lower bounds.
\end{proof}

\begin{remark}[Computing the enclosure]
\label{rem:finite-gain-computability}
The finite gain is the maximum singular value over a compact time interval. The defect uses finite solutions and coupling inner products, such as $\norm{Q_J^AC^*y}_H^2=\norm{C^*y}_H^2-\norm{R_JC^*y}_H^2$. Rigorous evaluation must bound the matrix exponential, the time maximum and these inner products. In the tail, use an upper bound for $\eta_{N'}$ and a positive lower bound for $\lambda_{N'+1}$, since $(1-e^{-\lambda\Tend})/\lambda=\int_0^{\Tend}e^{-\lambda s}\,\diff s$ decreases with $\lambda>0$. Floating-point values alone do not provide these guarantees.
\end{remark}

\begin{example}[A local coupling with explicit residual bounds]
\label{ex:finite-gain-local}
Let $\Omega=(0,\pi)$, $A=L=-\partial_{xx}$ with Dirichlet conditions, $Cw=2\cos(x)w$ and $\phi_j(x)=\sqrt{2/\pi}\sin(jx)$. We take $R_m=P_m$, $N<N'\le m$, and write $a_k(t),b_k(t)$ for the $m$th sine coefficients of $x_{m,m}^{(k)},y_{m,m}^{(k)}$. Because $C\phi_j=\phi_{j-1}+\phi_{j+1}$, with $\phi_0=0$, $\mathfrak r_{m,m}(t;\phi_k)=(b_k(t)\phi_{m+1},-a_k(t)\phi_{m+1})^{\top}$. Iterating the variation-of-constants formula, each coupling factor has norm at most two and changes the spatial index by at most one. The terms of order below $m-k$ therefore vanish in mode $m$, giving
\begin{align*}
 (|a_k(t)|^2+|b_k(t)|^2)^{1/2}
 \le\sum_{n=m-k}^{\infty}\frac{(2t)^n}{n!}
 \le e^{2\Tend}\frac{(2\Tend)^{m-k}}{(m-k)!}.
\end{align*}
Consequently,
\begin{align}\label{eq:finite-gain-local-explicit}
 \Delta_{m,m;N,N'}\le\Tend\sqrt{N'-N}
 \sum_{n=m-N'}^{\infty}\frac{(2\Tend)^n}{n!},
 \quad \overline\Gamma_{N',\Tend}\le\frac{2}{(N'+1)^2}.
\end{align}
For $\Tend=1$, $N=1$, $N'=10$, $m=26$, these give $\Delta_{26,26;1,10}\le3e^2 2^{16}/16!<7\times10^{-8}$ and $\overline\Gamma_{10,1}\le2/121$.

In sine coordinates ordered as $(x_1,\ldots,x_m,y_1,\ldots,y_m)^{\top}$, the finite generator is $G_m=\left(\begin{smallmatrix}-D_m&-B_m\\B_m&-D_m\end{smallmatrix}\right)$, where $D_m=\operatorname{diag}(1^2,\ldots,m^2)$ and $(B_m)_{ij}=1$ when $|i-j|=1$, and zero otherwise. Let $E:\R^9\to\R^{52}$ insert its coordinates into hidden modes $2,\ldots,10$, and let $\pi_u$ select the first 26 coordinates. We set $F(t)=\pi_u e^{tG_{26}}E$. All coordinate spaces have their Euclidean norms, so $\gamma_{26,26;1,10}=\max_{0\le t\le1} \norm{F(t)}_{\mathcal L(\R^9,\R^{26})}$. For a real matrix $Z=(z_{ij})$, we write $\norm{Z}_{\mathrm F}=(\sum_{i,j}|z_{ij}|^2)^{1/2}$ for its Frobenius norm.

Here is a fixed-point computation with explicit error bounds. We take $h=1/2000$, $D=2^{90}$ and $Z_0=DE$, and at step $i$ form integer matrices
\begin{align*}
 T_0=Z_i,\quad
 T_r=\operatorname{rnd}\bigl(G_{26}T_{r-1}/(2000r)\bigr)
 \ (1\le r\le12),\quad Z_{i+1}=\sum_{r=0}^{12}T_r,
\end{align*}
where $\operatorname{rnd}$ rounds each entry to the nearest integer, with ties towards $+\infty$. Exact integer comparisons give $\norm{\pi_u Z_i/D}_{\mathrm F}<0.158$ for $0\le i\le2000$ and $\norm{\pi_u Z_{780}e_1/D}_{\R^{26}}>0.1535$, where $e_1=(1,0,\ldots,0)^{\top}\in\R^9$. Because $\norm{G_{26}}_{\mathcal L(\R^{52})}\le678$ and $e^{678h}<3/2$, a Taylor step has remainder $r_*:=(3/2)(678h)^{13}/13!<3\times10^{-16}$. Each of its twelve roundings has matrix norm at most $\sqrt{52\cdot9}/(2D)$, and propagation through later stages gives a total rounding bound $\rho_*=18\sqrt{468}/(2D)$ per step. Contractivity of the exact evolution then bounds the grid error by
\begin{align*}
 \max_{0\le i\le2000}
 \norm{Z_i/D-e^{ihG_{26}}E}_{\mathcal L(\R^9,\R^{52})}
 \le2000(r_*+\rho_*)e^{2000r_*}<10^{-11}.
\end{align*}
The corresponding Frobenius error is below $3\times10^{-11}$.

For the piecewise affine time interpolant $I_h^tF$, contractivity and commutation with $G_{26}$ give
\begin{align*}
 \max_{0\le t\le1}\norm{F(t)-I_h^tF(t)}_{\mathrm F}
 \le\frac{h^2}{8}\max_{0\le t\le1}\norm{F''(t)}_{\mathrm F}
 \le\frac{h^2}{8}\norm{G_{26}^2E}_{\mathrm F}.
\end{align*}
An exact computation gives $\norm{G_{26}^2E}_{\mathrm F}^2=167838875<13000^2$. Convexity bounds the interpolant by its largest nodal norm, so
\begin{align*}
 0.1535-10^{-11}<\gamma_{26,26;1,10}
 <0.158+\frac{13000}{8\cdot2000^2}+3\times10^{-11}.
\end{align*}
Substitution in \eqref{eq:finite-gain-enclosure}, with rational comparisons, therefore proves
\begin{align}\label{eq:finite-gain-local-verified}
 0.15349<\Gamma_{1,\Tend}<0.1593\quad(\Tend=1).
\end{align}
The a priori estimate \eqref{eq:information-gain-upper} gives only $\Gamma_{1,1}\le(1-e^{-4})/2<1/2$. Thus, the optimal worst-case influence of an unknown hidden initial state of norm at most $R$ is below $0.1593R$, whereas no prediction from the known data can guarantee a uniform error below $0.15349R$.
\end{example}


\section{Two-sided decisions for known hidden data}
\label{sec:two-sided}
We assume that $v_0$ is fully known. The upper bound of \cref{thm:resolved-certificate} can certify accuracy but cannot reject a tolerance. Throughout this section, we write $\mathcal X:=\mathcal C([0,\Tend];H)$ with its maximum norm.

\begin{proposition}[Spectral data do not decide the tolerance]
\label{prop:spectral-only-impossibility}
There are coupled Dirichlet diffusion systems with the same $L,C,P_N$, $\lambda_{N+1},\eta_N$, initial and forcing data, and reduced trajectories, such that a fixed tolerance is met in one and violated in the other.
\end{proposition}

\begin{proof}
On $(0,\pi)$, we take $L=-\partial_{xx}$, $C=I$ and $A_a=-\partial_{xx}+aI$, $a\ge0$, with Dirichlet conditions. Let $u_0=0$, $v_0=\phi_j$, $j=N+1$, and $f=g=0$, where $\phi_j=\sqrt{2/\pi}\sin(jx)$. For any $a$, the reduction gives $u_N=v_N=0$, while $\lambda_{N+1}=j^2$ and $\eta_N=1$. For $a=0$, $u(t)=-e^{-j^2t}\sin(t)\phi_j$. For general $a$, the energy identity bounds $\norm{v(t)}_H$ by one; both fields remain in mode $j$, so variation of constants gives $\norm{u}_{\mathcal X}\le1/(j^2+a)$. We choose $N=1$, $\Tend=1$ and $\varepsilon=1/20$. At $a=0$, $\norm{u(1/4)}_H=e^{-1}\sin(1/4)>1/20$; at $a=20$, $\norm{u}_{\mathcal X}\le1/24<1/20$.
\end{proof}

\begin{remark}\label{rem:information-boundary-scope}
Spectral damping and coupling norms provide sufficient bounds, but \cref{prop:spectral-only-impossibility} rules out a necessary and sufficient tolerance test based on them alone. The following auxiliary problem also uses the retained generator $A$.
\end{remark}

\subsection{A lower bound from an auxiliary problem}
\label{subsec:two-sided-decision}
In the setting of \cref{thm:resolved-certificate}, let $e_N=u-u_N$, $d_N=C^*a_N$ and $K(t)=C^*e^{-tL}C$. The transmitted omitted response $d_N$ belongs to $L^2(0,\Tend;H)$, since $a_N\in\W(0,\Tend)$. We define the auxiliary response
\begin{align}\label{eq:decision-direct-response}
 \mathfrak F_N(t):=-\int_0^t e^{-(t-s)A}d_N(s)\,\diff s,
 \quad 0\le t\le\Tend.
\end{align}
It lies in $\W(0,\Tend)$ and solves $\mathfrak F_N'+A\mathfrak F_N=-d_N$ in $V^*$ almost everywhere, with zero initial value. Its computation uses $u_N$ and the data, not $u$.

We put $\mu_A=\inf\sigma(A)>0$, $\mu_L=\inf\sigma(L)>0$ and choose a nonnegative $k\in L^1(0,\Tend)$ with $\norm{K(t)}_{\mathcal L(H)}\le k(t)$ almost everywhere. We define
\begin{align}\label{eq:decision-feedback-gain}
 \beta_{\Tend}:=\int_0^{\Tend}k(\tau)
 \frac{1-e^{-\mu_A(\Tend-\tau)}}{\mu_A}\,\diff\tau.
\end{align}
For the admissible choice $k(\tau)=\norm{C}_{\mathcal L(H)}^2e^{-\mu_L\tau}$,
\begin{align}\label{eq:decision-explicit-gain}
 \beta_{\Tend}
 =\frac{\norm{C}_{\mathcal L(H)}^2}{\mu_A}
 \int_0^{\Tend}e^{-\mu_L\tau}
 (1-e^{-\mu_A(\Tend-\tau)})\,\diff\tau
 \le\frac{\norm{C}_{\mathcal L(H)}^2}{\mu_A\mu_L}.
\end{align}
The equality in \eqref{eq:decision-explicit-gain} concerns this choice of $k$; an arbitrary majorant need not give the same or a smaller bound.

\begin{theorem}[Two-sided bound and three-way test]
\label{thm:two-sided-retained-decision}
Under the hypotheses of \cref{thm:resolved-certificate}, let the quantities
above be given, and take $U_N\ge0$ such that
\begin{align}\label{eq:decision-known-upper}
 \norm{e_N}_{\mathcal X}\le U_N.
\end{align}
For example, we take $U_N=\norm{d_N}_{L^1(0,\Tend;H)}$ or $U_N=\mathcal E_N$ from \eqref{eq:certificate-explicit}. We set $D_N=\norm{\mathfrak F_N}_{\mathcal X}$ and
\begin{align}
 \underline E_N^{\mathrm{tr}}
 &:=\max\left\{0,\frac{D_N}{1+\beta_{\Tend}},D_N-\beta_{\Tend}U_N\right\},
 \label{eq:decision-lower-bound}\\
 U_N^\sharp&:=\min\{U_N,D_N+\beta_{\Tend}U_N\}.
 \label{eq:decision-sharp-upper}
\end{align}
If $\beta_{\Tend}<1$, also include $D_N/(1-\beta_{\Tend})$ in the minimum. Then,
\begin{align}\label{eq:decision-two-sided}
 \underline E_N^{\mathrm{tr}}\le\norm{u-u_N}_{\mathcal X}\le U_N^\sharp.
\end{align}
For $\varepsilon>0$, accept if $U_N^\sharp\le\varepsilon$, reject if $\underline E_N^{\mathrm{tr}}>\varepsilon$, and otherwise leave the test undecided. These decisions concern the chosen model and tolerance, not dynamical stability. Furthermore, $u=u_N$ if and only if $\mathfrak F_N=0$.
\end{theorem}

\begin{proof}
Let $z_N=v-v_N-a_N$. By the proof of \cref{thm:resolved-certificate}, $z_N'+Lz_N=Ce_N$, $z_N(0)=0$. Eliminating $z_N$ gives
\begin{align}\label{eq:decision-error-volterra}
 e_N'+Ae_N+K*e_N=-d_N,\quad e_N(0)=0.
\end{align}
All convolutions are well-defined in $H$. Variation of constants yields
\begin{align}\label{eq:decision-fixed-point}
 e_N=\mathfrak F_N-\mathcal B_{\Tend}e_N,\quad
 (\mathcal B_{\Tend}w)(t)
 :=\int_0^t e^{-(t-s)A}(K*w)(s)\,\diff s.
\end{align}
For $w\in\mathcal X$, Tonelli's theorem and the semigroup bound give
\begin{align*}
 \norm{(\mathcal B_{\Tend}w)(t)}_H
 &\le\norm{w}_{\mathcal X}\int_0^t k(\tau)
 \frac{1-e^{-\mu_A(t-\tau)}}{\mu_A}\,\diff\tau
 \le\beta_{\Tend}\norm{w}_{\mathcal X}.
\end{align*}
Thus,
\begin{align}\label{eq:decision-feedback-difference}
 \norm{e_N-\mathfrak F_N}_{\mathcal X}
 \le\beta_{\Tend}\norm{e_N}_{\mathcal X}\le\beta_{\Tend}U_N.
\end{align}
The triangle inequality also gives $D_N\le(1+\beta_{\Tend})\norm{e_N}_{\mathcal X}$ and, if $\beta_{\Tend}<1$, $(1-\beta_{\Tend})\norm{e_N}_{\mathcal X}\le D_N$. These inequalities prove the endpoints. If $e_N=0$, then \eqref{eq:decision-fixed-point} gives $\mathfrak F_N=0$. Conversely, $\mathfrak F_N=0$ implies $d_N=0$ almost everywhere by its variational equation. The zero-initial-state passive system for $(e_N,z_N)^{\top}$ then gives $e_N=0$, without requiring $\beta_{\Tend}<1$.
\end{proof}

\begin{remark}\label{rem:decision-limitations}
For the exact auxiliary response and finite $\beta_{\Tend}$, a positive error implies $D_N>0$ and hence $\underline E_N^{\mathrm{tr}}>0$. Nevertheless, the lower bound may be too small to reject, or the upper bound too large to accept. The parameter $\beta_{\Tend}$ bounds feedback; it is not a stability threshold. The basic endpoints remain valid for $\beta_{\Tend}\ge1$.
\end{remark}

\begin{corollary}[Approximate auxiliary solution]
\label{cor:three-way-with-auxiliary-error}
Under the hypotheses of \cref{thm:two-sided-retained-decision}, suppose $\mathfrak F_N^{\rm app}\in\mathcal X$ and $\delta_N\ge0$ satisfy $\norm{\mathfrak F_N-\mathfrak F_N^{\rm app}}_{\mathcal X}\le\delta_N$. We set
\begin{align*}
 D_N^-&:=\max\{0,\norm{\mathfrak F_N^{\rm app}}_{\mathcal X}-\delta_N\},
 \quad D_N^+:=\norm{\mathfrak F_N^{\rm app}}_{\mathcal X}+\delta_N,\\
 \underline E_N&:=\max\{0,D_N^-/(1+\beta_{\Tend}),D_N^--\beta_{\Tend}U_N\},
 \quad \overline E_N:=\min\{U_N,D_N^++\beta_{\Tend}U_N\}.
\end{align*}
If $\beta_{\Tend}<1$, also include $D_N^+/(1-\beta_{\Tend})$ in the minimum. Then,
\begin{align}\label{eq:decision-certified-interval}
 \underline E_N\le\norm{u-u_N}_{\mathcal X}\le\overline E_N.
\end{align}
Accept if $\overline E_N\le\varepsilon$, reject if $\underline E_N>\varepsilon$, and otherwise leave the test undecided. All scalar comparisons, including the trajectory norm, require rigorous bounds for a numerical decision. Here the lower endpoint can vanish for a positive error because of the auxiliary error allowance.
\end{corollary}
\begin{proof}
The reverse triangle inequality gives $D_N^-\le D_N\le D_N^+$. We use monotonicity in $D_N$ in \cref{thm:two-sided-retained-decision}.
\end{proof}

\subsection{Discretisation}
\label{subsec:fully-discrete-certification}
We retain the spectral hidden space $H_N$ and distinguish reduction error from numerical error. Fix $N\ge1$, we put $\eta_N=\norm{Q_NC}_{\mathcal L(H)}$, and let $(u_N,v_N)^{\top}$ be the exact reduced solution. Take continuous time reconstructions $U_h\in\W(0,\Tend)$, $V_h\in H^1(0,\Tend;H_N)$ and we define
\begin{align}
 R_u&:=f-U_h'-AU_h-C^*V_h\in L^2(0,\Tend;V^*),
 \label{eq:discrete-residual-u}\\
 R_v&:=P_Ng-V_h'-LV_h+P_NCU_h\in L^2(0,\Tend;H_N).
 \label{eq:discrete-residual-v}
\end{align}
The operator $A$ is understood variationally, and $LV_h\in H_N$. Conforming piecewise affine spatial approximations are therefore allowed. Writing $e_h^u=u_N-U_h$, we require a bound $\Delta_N\ge0$ such that
\begin{align}\label{eq:discrete-reduced-certified-budget}
 \norm{e_h^u}_{\mathcal X}\le\Delta_N.
\end{align}
One choice is an upper bound for the nonnegative square root defined by
\begin{align}\label{eq:discrete-reduced-budget}
 \begin{aligned}
 (\Delta_N^{\mathrm{res}})^2
 :={}&\norm{u_0-U_h(0)}_H^2+\norm{P_Nv_0-V_h(0)}_H^2\\
 &+\alpha_A^{-1}\norm{R_u}_{L^2(0,\Tend;V^*)}^2
 +\lambda_1^{-1}\norm{R_v}_{L^2(0,\Tend;H)}^2.
 \end{aligned}
\end{align}
An alternative rigorous a priori or a posteriori bound may be used.

The omitted response driven by $U_h$ is
\begin{align}
 \widetilde a_N(t)&:=e^{-tL}Q_Nv_0
 +\int_0^t e^{-(t-s)L}Q_N(CU_h(s)+g(s))\,\diff s,
 \label{eq:discrete-omitted-response}\\
 \widetilde d_N&:=C^*\widetilde a_N.
\end{align}
We choose $d_h\in L^2(0,\Tend;H)$ with a bound $\xi_N\ge0$ satisfying
\begin{align}\label{eq:discrete-response-budget}
 \norm{\widetilde d_N-d_h}_{L^1(0,\Tend;H)}\le\xi_N.
\end{align}
Let $F$ solve $F'+AF=-d_h$, $F(0)=0$. For its reconstruction $F_h\in\W(0,\Tend)$, we take $\delta_F\ge0$ such that
\begin{align}\label{eq:discrete-auxiliary-certified-budget}
 \norm{F-F_h}_{\mathcal X}\le\delta_F.
\end{align}
For example, one may bound the nonnegative square root given by
\begin{align}\label{eq:discrete-auxiliary-budget}
 \begin{aligned}
 R_F&:=-d_h-F_h'-AF_h\in L^2(0,\Tend;V^*),\\
 (\delta_F^{\mathrm{res}})^2
 &:=\norm{F_h(0)}_H^2+\alpha_A^{-1}\norm{R_F}_{L^2(0,\Tend;V^*)}^2.
 \end{aligned}
\end{align}

\begin{theorem}[Decisions with discretisation errors]
\label{thm:discrete-three-way}
Assume \cref{ass:forms} and \eqref{eq:data-full}, and let the preceding reconstructions and bounds be given. We take $\overline d_h,D_h^-,D_h^+\ge0$ such that $\norm{d_h}_{L^1(0,\Tend;H)}\le\overline d_h$ and $D_h^-\le\norm{F_h}_{\mathcal X}\le D_h^+$. Suppose $\overline\beta\ge0$ satisfies
\begin{align}\label{eq:discrete-restricted-feedback}
 \norm{\mathcal B_{\Tend}e_N}_{\mathcal X}
 \le\overline\beta\norm{e_N}_{\mathcal X},\quad e_N=u-u_N.
\end{align}
Any $\overline\beta\ge\beta_{\Tend}$ is admissible. A smaller value may be used on an invariant subspace containing $e_N$, if this restriction is proved from the equations and data. We define
\begin{align}
 \kappa_N&:=\frac{\eta_N^2\Tend}{\lambda_{N+1}},
 \quad \vartheta_N :=\kappa_N\Delta_N+\xi_N+\delta_F,
 \label{eq:discrete-total-budget}\\
 \overline d_N&:=\overline d_h+\kappa_N\Delta_N+\xi_N.
 \label{eq:discrete-forcing-upper}
\end{align}
A rigorous upper bound for $\kappa_N$ may replace its exact value throughout. For $E_N=\norm{u-u_N}_{\mathcal X}$,
\begin{align}\label{eq:discrete-decision-enclosure}
 \begin{aligned}
 \underline E_N
 &:=\max\left\{0,\frac{D_h^--\vartheta_N}{1+\overline\beta},
 D_h^--\vartheta_N-\overline\beta\,\overline d_N\right\}\le E_N,\\
 E_N&\le\overline E_N
 :=\min\{\overline d_N,D_h^++\vartheta_N+\overline\beta\,\overline d_N\}.
 \end{aligned}
\end{align}
If $\overline\beta<1$, also include $(D_h^++\vartheta_N)/(1-\overline\beta)$ in the minimum. Furthermore,
\begin{align}\label{eq:discrete-total-retained-error}
 \max\{0,\underline E_N-\Delta_N\}
 \le\norm{u-U_h}_{\mathcal X}\le\overline E_N+\Delta_N.
\end{align}
Both intervals give three-way tolerance tests as in \cref{thm:two-sided-retained-decision}.
\end{theorem}

\begin{proof}
We set $r=u_N-U_h$, $s=v_N-V_h$. Cancellation of the coupling in the product energy and Young's inequality give, almost everywhere,
\begin{align*}
 \frac{\diff}{\diff t}(\norm{r}_H^2+\norm{s}_H^2)
 +\alpha_A\norm{r}_V^2+\lambda_1\norm{s}_H^2
 \le\alpha_A^{-1}\norm{R_u}_{V^*}^2+\lambda_1^{-1}\norm{R_v}_H^2.
\end{align*}
Integration proves the residual estimate \eqref{eq:discrete-reduced-budget}. The same argument for $F-F_h$ proves \eqref{eq:discrete-auxiliary-budget}. More generally we use
\begin{align}\label{eq:discrete-reduced-reliability}
 \norm{u_N-U_h}_{\mathcal X}\le\Delta_N.
\end{align}
The spectral tail estimate and Tonelli's theorem give
\begin{align*}
 \norm{d_N-\widetilde d_N}_{L^1(0,\Tend;H)}
 &\le\frac{\eta_N^2}{\lambda_{N+1}}\int_0^{\Tend}
 (1-e^{-\lambda_{N+1}(\Tend-t)})\norm{u_N(t)-U_h(t)}_H\,\diff t
 \le\kappa_N\Delta_N.
\end{align*}
Consequently,
\begin{align}\label{eq:discrete-transmission-bound}
 \norm{d_N-d_h}_{L^1(0,\Tend;H)}\le\kappa_N\Delta_N+\xi_N,
 \quad d_N=C^*a_N.
\end{align}
By \eqref{eq:certificate-basic}, 
\begin{align*}
E_N\le\norm{d_N}_{L^1(0,\Tend;H)} \le\overline d_N
\end{align*}
Contractivity also gives $\norm{\mathfrak F_N-F}_{\mathcal X} \le\norm{d_N-d_h}_{L^1(0,\Tend;H)}$, hence
\begin{align*}
 \norm{\mathfrak F_N-F_h}_{\mathcal X}\le\vartheta_N,\quad
 \max\{0,D_h^--\vartheta_N\}\le D_N\le D_h^++\vartheta_N.
\end{align*}
The proof of \cref{thm:two-sided-retained-decision} uses the feedback bound only for $e_N$. Applying it with \eqref{eq:discrete-restricted-feedback} proves \eqref{eq:discrete-decision-enclosure}; the triangle inequalities and \eqref{eq:discrete-reduced-reliability} prove \eqref{eq:discrete-total-retained-error}.
\end{proof}

\begin{remark}[Evaluating the norms]
\label{rem:finite-nodal-norm-certification}
Let $0=t_0<\cdots<t_m=\Tend$ and suppose $X_h\in\mathcal X$ is affine on each time interval. Convexity gives
\begin{align}\label{eq:nodal-C-norm}
 \norm{X_h}_{\mathcal X}=\max_{0\le n\le m}\norm{X_h(t_n)}_H.
\end{align}
For a finite-element basis $\{\psi_i\}_{i=1}^r$, we write $X_h(t_n)=\sum_i c_{n,i}\psi_i$, $c_n=(c_{n,1},\ldots,c_{n,r})^{\top}$, and $M=((\psi_i,\psi_j)_H)_{i,j=1}^r$. Then,
\begin{align}\label{eq:nodal-mass-norm}
 \norm{X_h(t_n)}_H^2=c_n^{\top}Mc_n.
\end{align}
Verified bounds for these quadratic forms, including the mass entries, give $D_h^\pm$ when $X_h=F_h$. A single node suffices for a lower bound. If $d_h$ is also affine on each interval, then
\begin{align}\label{eq:nodal-L1-bound}
 \norm{d_h}_{L^1(0,\Tend;H)}\le
 \sum_{n=1}^m\frac{t_n-t_{n-1}}2
 (\norm{d_h(t_{n-1})}_H+\norm{d_h(t_n)}_H).
\end{align}
For other forcing reconstructions the integral needs a separate bound. These norm evaluations do not supply $\Delta_N,\xi_N,\delta_F$.
\end{remark}

\begin{remark}[Spatial and temporal residuals]
\label{rem:discrete-flux}
For a bounded symmetric uniformly positive definite diffusion matrix field $\mathsf A$, let $a_A(w,\varphi)=(\mathsf A\nabla w,\nabla\varphi)_{L^2(\Omega)^d}$. For the following $L^2$ flux bound, assume also $U_h'\in L^2(0,\Tend;H)$, as holds for a finite-element reconstruction that is piecewise affine in time. We define $H(\operatorname{div};\Omega)$ as the vector fields in $L^2(\Omega)^d$ whose distributional divergence belongs to $L^2(\Omega)$. For any reconstructed flux $\sigma_h\in L^2(0,\Tend;H(\operatorname{div};\Omega))$, integration by parts against $V=H_0^1(\Omega)$ gives
\begin{align*}
 \norm{R_u(t)}_{V^*}
 \le C_P\norm{f-U_h'-C^*V_h+\operatorname{div}\sigma_h}_{L^2(\Omega)}
 +\norm{\sigma_h-\mathsf A\nabla U_h}_{L^2(\Omega)^d},
\end{align*}
where $\norm{w}_H\le C_P\norm{w}_V$ for $w\in V$. The same construction applies to $R_F$ when $F_h'\in L^2(0,\Tend;H)$. These estimates require rigorous bounds for data and quadrature errors; a finite-element computation alone does not supply them. The next example instead uses analytic a priori bounds. Non-invariant hidden trial spaces would require additional projection-defect estimates and are outside the present theorem.
\end{remark}

\subsection{A coupling that excites infinitely many modes}
\label{subsec:local-multiplication-cascade}

\begin{example}[Mode mixing by a multiplication coupling]
\label{ex:local-cascade-three-way}
Let $H=L^2(0,\pi)$, $A=L=-\partial_{xx}$ with Dirichlet conditions, and $\phi_j(x)=\sqrt{2/\pi}\sin(jx)$, $j\ge1$. We take
\begin{align}\label{eq:local-coupling}
 (Cw)(x)=2\cos(x)w(x),\quad w\in H.
\end{align}
Then, $C=C^*$, $\norm{C}_{\mathcal L(H)}=2$, and
\begin{align}\label{eq:local-neighbour-mixing}
 C\phi_1=\phi_2,\quad C\phi_j=\phi_{j-1}+\phi_{j+1}\quad(j\ge2).
\end{align}
The orthonormal sequence $\phi_j$ is weakly null, whereas $\norm{C\phi_j}_H=\sqrt2$ for $j\ge2$, so $C$ is not compact. Let $u_0=\phi_1$, $v_0=\rho\phi_2$, $\rho\ge0$, $f=g=0$, $\Tend=1$ and $H_1=\operatorname{span}\{\phi_1\}$. Since $P_1C\phi_1=P_1v_0=0$,
\begin{align}\label{eq:local-exact-reduction}
 u_1(t)=e^{-t}\phi_1,\quad v_1(t)=0.
\end{align}
The full trajectory cannot stay in a fixed finite span of sine modes. Its generator on $D(A)\times D(L)$ is $\mathcal G U=(-Ax-C^*y,Cx-Ly)^{\top}$, $U=(x,y)^{\top}$. It maps finite sine sums to finite sums, so $U_0=(\phi_1,\rho\phi_2)^{\top}$ lies in every $D(\mathcal G^k)$ and $U^{(k)}(0)=\mathcal G^kU_0$. Its largest mode is $k+1$ when $\rho=0$ and $k+2$ when $\rho>0$: each new largest index is reached by a unique chain of coupling factors $\pm1$, with nonzero initial coefficient. These derivatives cannot all belong to a fixed finite span.

Nevertheless, the omitted and transmitted responses are explicit:
\begin{align}\label{eq:local-transmitted-response}
 \begin{aligned}
 a_1(t)&=q_\rho(t)\phi_2,\quad
 d_1(t)=q_\rho(t)(\phi_1+\phi_3),\\
 q_\rho(t)&=\rho e^{-4t}+(e^{-t}-e^{-4t})/3\ge0.
 \end{aligned}
\end{align}
We set
\begin{align}\label{eq:local-q-integral}
 I_\rho:=\int_0^1 q_\rho(t)\,\diff t
 =\frac{\rho(1-e^{-4})}{4}
 +\frac13\left[1-e^{-1}-\frac{1-e^{-4}}4\right].
\end{align}
The auxiliary response is
\begin{align}\label{eq:local-auxiliary-formula}
 \mathfrak F_1(t)=-G_{1,\rho}(t)\phi_1-G_{3,\rho}(t)\phi_3,
\end{align}
where
\begin{align}
 G_{1,\rho}(t)&=\frac\rho3(e^{-t}-e^{-4t})
 +\frac{te^{-t}}3-\frac{e^{-t}-e^{-4t}}9,
 \label{eq:local-g1}\\
 G_{3,\rho}(t)&=\frac\rho5(e^{-4t}-e^{-9t})
 +\frac13\left[\frac{e^{-t}-e^{-9t}}8-\frac{e^{-4t}-e^{-9t}}5\right].
 \label{eq:local-g3}
\end{align}
Both coefficients equal $\int_0^t e^{-j^2(t-s)}q_\rho(s)\,\diff s$ for $j=1,3$, and are nonnegative.

Let $H_{\rm odd}$ and $H_{\rm even}$ be the closed spans of the odd- and even-index sine modes. The coupled equations and these data preserve $H_{\rm odd}\times H_{\rm even}$, so $e_1\in\mathcal C([0,1];H_{\rm odd})$. Because $C$ interchanges the two spaces and $L\ge4I$ on $H_{\rm even}$,
\begin{align}\label{eq:local-feedback-gain}
 \norm{K(t)|_{H_{\rm odd}}}_{\mathcal L(H_{\rm odd})}\le4e^{-4t},\quad
 \beta_{1,\rm odd}:=
 4\left[\frac{1-e^{-4}}4-\frac{e^{-1}-e^{-4}}3\right]<\frac{13}{25}.
\end{align}
This bounds feedback on $H_{\rm odd}$; it is not a claim about the norm on all of $H$. Thus, \eqref{eq:discrete-restricted-feedback} holds for the actual error with $\overline\beta=13/25$, and
\begin{align}\label{eq:local-continuous-enclosure}
 \frac{G_{1,\rho}(1)}{1+13/25}
 \le\norm{u-u_1}_{\mathcal X}\le\sqrt2\,I_\rho.
\end{align}
No full trajectory is used in either endpoint.

For a numerical approximation, divide $(0,\pi)$ into $J=256$ equal intervals of length $h=\pi/J$. Use continuous piecewise affine elements with zero boundary values and the consistent $L^2$ inner product. Let $I_h$ be nodal interpolation and put $\chi_j=I_h\phi_j$, $j=1,3$. These are orthogonal discrete sine eigenvectors for the consistent mass matrix, with
\begin{align}\label{eq:local-fe-eigenvalues}
 \mu_j=\frac{6(1-\cos(jh))}{h^2(2+\cos(jh))},\quad
 \nu_j:=\frac{(\phi_j,\chi_j)_H}{\norm{\chi_j}_H^2}=\frac{\mu_j}{j^2}.
\end{align}
Indeed, with $z=jh$, direct integration gives $\norm{\chi_j}_H^2=(2+\cos z)/3$ and $(\phi_j,\chi_j)_H=2(1-\cos z)/(j^2h^2)$; the stiffness and mass eigenvalues give $\mu_j$. The formula for $\nu_j$ does not apply to a lumped mass matrix. We take $m=4096$, $\tau=1/m$, $t_n=n\tau$. The backward Euler reduced field, reconstructed affinely in time, has nodes
\begin{align}\label{eq:local-reduced-be}
 U_{h,\tau}(t_n)=(1+\tau\mu_1)^{-n}\chi_1.
\end{align}
Its hidden coordinate stays zero because $(C\chi_1,\phi_1)_H=0$ by reflection about $\pi/2$. Spatial interpolation and scalar backward Euler stability give
\begin{align}\label{eq:local-reduced-budget}
 \norm{u_1-U_{h,\tau}}_{\mathcal X}
 \le\frac1{J^2}+\frac{h^2}{5}
 +\left(\frac1{2m}+\frac1{8m^2}\right)\left(\frac{1001}{1000}\right)^2
 <\frac1{400}=: \Delta_1.
\end{align}

We use the exact known forcing $d_h=d_1$. Its finite-element backward Euler solution $F_{h,\tau}$ is the affine time reconstruction of
\begin{align}\label{eq:local-auxiliary-be}
 \begin{aligned}
 F_{h,\tau}(t_n)&=y_{1,n}\chi_1+y_{3,n}\chi_3,\quad y_{j,0}=0,\\
 \frac{y_{j,n}-y_{j,n-1}}\tau+\mu_jy_{j,n}
 &=-\nu_jq_\rho(t_n),\quad j=1,3.
 \end{aligned}
\end{align}
The full finite-element load has only these two discrete sine components, so no other coordinates occur. Uniformly for $0\le\rho\le3$,
\begin{align}\label{eq:local-auxiliary-budget}
 \norm{\mathfrak F_1-F_{h,\tau}}_{\mathcal X}<\frac1{64}=: \delta_F.
\end{align}
To justify both error budgets, note that $h<1/80$ and, for $j=1,3$,
\begin{align*}
 \begin{gathered}
 1\le\nu_j\le1+(jh)^2/5,\quad
 \mu_1<1001/1000,\quad\mu_3<901/100,\\
 \norm{\phi_j-\chi_j}_H\le j^2/J^2,\quad\norm{\chi_j}_H\le1.
 \end{gathered}
\end{align*}
For $0<z=jh<3/80$, the alternating bounds for $\cos z$ imply $1\le6(1-\cos z)/(z^2(2+\cos z))\le(1-z^2/6)^{-1}\le1+z^2/5$. On each element the interpolation error vanishes at both endpoints and its derivative has mean zero. Two applications of the one-dimensional Poincar\'e inequality, for zero endpoint values and for zero mean, each with constant $h/\pi$, give the interpolation bound. The formula for $\norm{\chi_j}_H$ gives its upper bound.

Writing $Q_0=10/3$ and $Q_1=41/3$, we have $|q_\rho|\le Q_0$ and $|q_\rho'|\le Q_1$ on $[0,1]$. The semidiscrete coefficient $y_j'+\mu_j y_j=-\nu_jq_\rho$, $y_j(0)=0$, satisfies $\norm{y_j+G_{j,\rho}}_{\mathcal C([0,1];\R)}\le2Q_0(jh)^2/5$. Also $|y_j|\le Q_0/j^2$, $|y_j'|\le2\nu_jQ_0$ and $|y_j''|\le\nu_j(Q_1+2\mu_jQ_0)$, hence $|y_1''|<21$ and $|y_3''|<75$. For a scalar equation with positive decay, backward Euler nodal error is at most $(\tau/2)\max|y''|$ on $[0,1]$; affine reconstruction adds $(\tau^2/8)\max|y''|$. Applied to $e^{-\mu_1t}$ this proves \eqref{eq:local-reduced-budget}. Applied to $y_1,y_3$, together with $G_{j,\rho}\le Q_0/j^2$, it gives
\begin{align}\label{eq:local-auxiliary-budget-proof}
 \begin{aligned}
 \norm{\mathfrak F_1-F_{h,\tau}}_{\mathcal X}
 &\le\frac{2(10/3)(1^2+3^2)}{5\cdot80^2}
 +\frac{2(10/3)}{J^2}\\
 &\quad +(21+75)\left(\frac1{2m}+\frac1{8m^2}\right)<\frac1{64}.
 \end{aligned}
\end{align}
The three terms bound the coefficient error, spatial interpolation and time discretisation, respectively.

Although $d_h=d_1$, the forcing $\widetilde d_1$ driven by $U_{h,\tau}$ differs from $d_h$. The spectral comparison gives
\begin{align}\label{eq:local-transmission-budget}
 \norm{\widetilde d_1-d_h}_{L^1(0,1;H)}
 \le\frac{\norm{C}_{\mathcal L(H)}^2}{\lambda_2}
 \norm{U_{h,\tau}-u_1}_{\mathcal X}\le\frac1{400}=: \xi_1.
\end{align}
We take $\overline d_h=(17/12)I_\rho$, since $\sqrt2<17/12$. Instead of evaluating the finite nodal norms, use the valid but conservative brackets $D_h^-=\max\{0,G_{1,\rho}(1)-1/64\}$ and $D_h^+=\overline d_h+1/64$. With an upper bound $\kappa_1=1$, the budgets in \cref{thm:discrete-three-way} are $\overline\beta=13/25$, $\overline d_1=\overline d_h+1/200$ and $\vartheta_1=1/200+1/64$. Here, $\delta_F$ is used both in the nodal brackets and in $\vartheta_1$; this is conservative, not an additional source of error.

At $\varepsilon=1/5$, Taylor bounds for the exponentials and rational comparisons give the following decisions for $U_{h,\tau}$:
\begin{center}
\begin{tabular}{@{}lll@{}}
\toprule
Amplitude & Outcome & Reason\\
\midrule
$\rho=0$ & Accept & $\overline d_h+1/200+1/400<1/5$\\
$\rho=1$ & Undecided & the certified interval contains $1/5$\\
$\rho=3$ & Reject &
$\displaystyle\frac{G_{1,3}(1)-1/64-\vartheta_1}{1+13/25}
-\frac1{400}>\frac15$\\
\bottomrule
\end{tabular}
\end{center}
Thus, the test applies to a local non-compact coupling whose full dynamics cannot remain in any fixed finite set of sine modes. It uses the reduced field, an auxiliary finite-element diffusion problem and bounds for the omitted transmission, without computing the full trajectory.
\end{example}

\section{Conclusion}
For known hidden initial data, the discarded response, damped by diffusion and transmitted through the coupling, bounds the retained-trajectory error from above. An auxiliary diffusion problem gives a lower bound. With rigorous discretisation estimates, these bounds yield acceptance, rejection or an undecided outcome at a prescribed tolerance.

For missing initial data in a norm ball of radius $R$, the radius of information $R\Gamma_{N,\Tend}$ is the optimal worst-case error over all admissible states. Our main result encloses the time-uniform gain by a finite spectral evolution, its full-equation residual and a diffusion tail. The enclosure converges without compactness or smallness of the coupling. Strict comparisons of this radius with a tolerance have finite certificates provided the endpoint quantities admit arbitrarily accurate rigorous enclosures; nothing is claimed at equality. Together with the known-data bounds, this can certify a tolerance for all admissible states or show that some admissible state violates it.

The multiplication-coupling examples give explicit enclosures and finite-element error budgets for dynamics involving infinitely many modes. Nonlinear hidden dynamics, truncations not invariant under $L$, and a general-purpose verified implementation remain open.

\section*{Statements and declarations}

\subsection*{Funding}
The author declares that no funds, grants, or other support were received during the preparation of this manuscript.

\subsection*{Data availability}
No external datasets were used. The coefficients, initial data and discretisation parameters of the examples are specified in the manuscript. No additional numerical datasets or verification code are deposited or provided as supplementary material. The numerical inequalities are supported by the analytical error bounds stated in the text.

\subsection*{Competing interests}
The author declares no competing interests.

\subsection*{Use of generative AI}
The research questions, mathematical results, and conclusions are the author's own. The author used ChatGPT (OpenAI) for language editing, assistance with checking selected mathematical derivations and bibliographic information, and manuscript presentation. The author reviewed the AI-assisted work and takes full responsibility for the manuscript.



{\small
\begin{thebibliography}{99}

\bibitem{BeattieGugercinMehrmann2017}
Beattie, C., Gugercin, S., Mehrmann, V.: Model reduction for systems with inhomogeneous initial conditions. \emph{Systems Control Lett.} \textbf{99}, 99--106 (2017). \url{https://doi.org/10.1016/j.sysconle.2016.11.007}

\bibitem{BeckerHartmann2019}
Becker, S., Hartmann, C.: Infinite-dimensional bilinear and stochastic balanced truncation with explicit error bounds. \emph{Math. Control Signals Syst.} \textbf{31}, 1--37 (2019). \url{https://doi.org/10.1007/s00498-019-0234-8}

\bibitem{BreitenUnger2022}
Breiten, T., Unger, B.: Passivity preserving model reduction via spectral factorization. \emph{Automatica} \textbf{142}, 110368 (2022). \url{https://doi.org/10.1016/j.automatica.2022.110368}

\bibitem{CurtainZwart2020}
Curtain, R., Zwart, H.: \emph{Introduction to Infinite-Dimensional Systems Theory: A State-Space Approach}. Texts in Applied Mathematics. Springer, New York (2020). \url{https://doi.org/10.1007/978-1-0716-0590-5}

\bibitem{EngelNagel2000}
Engel, K.-J., Nagel, R.: \emph{One-Parameter Semigroups for Linear Evolution Equations}. Graduate Texts in Mathematics, vol.~194. Springer, New York (2000). \url{https://doi.org/10.1007/b97696}

\bibitem{FigotinSchenker2007}
Figotin, A., Schenker, J.H.: Hamiltonian treatment of time dispersive and dissipative media within the linear response theory. \emph{J. Comput. Appl. Math.} \textbf{204}, 199--208 (2007). \url{https://doi.org/10.1016/j.cam.2006.01.038}

\bibitem{GloverCurtainPartington1988}
Glover, K., Curtain, R.F., Partington, J.R.: Realisation and approximation
of linear infinite-dimensional systems with error bounds.
\emph{SIAM J. Control Optim.} \textbf{26}, 863--898 (1988).
\url{https://doi.org/10.1137/0326049}

\bibitem{GreplPatera2005}
Grepl, M.A., Patera, A.T.: A posteriori error bounds for reduced-basis approximations of parametrized parabolic partial differential equations. \emph{ESAIM Math. Model. Numer. Anal.} \textbf{39}, 157--181 (2005). \url{https://doi.org/10.1051/m2an:2005006}

\bibitem{GrubisicLazar2026}
Grubi\v{s}i\'c, L., Lazar, M.: Pollution free eigenvalue bounds for the Gramian operator. \emph{Syst. Control Lett.} \textbf{208}, 106308 (2026). \url{https://doi.org/10.1016/j.sysconle.2025.106308}

\bibitem{HeinkenschlossReisAntoulas2011}
Heinkenschloss, M., Reis, T., Antoulas, A.C.: Balanced truncation model reduction for systems with inhomogeneous initial conditions. \emph{Automatica} \textbf{47}(3), 559--564 (2011). \url{https://doi.org/10.1016/j.automatica.2010.12.002}

\bibitem{IshizakaCoercivity2026}
Ishizaka, H.: Coercivity structure of positive-type memory: exact gaps, critical horizons, and singular limits. Preprint, arXiv:2607.12482 (2026). \url{https://arxiv.org/abs/2607.12482}

\bibitem{IshizakaGraphSpace2026}
Ishizaka, H.: Graph-space well-posedness for diffusion equations with degenerate instantaneous diffusion. Preprint, arXiv:2607.12871 (2026). \url{https://arxiv.org/abs/2607.12871}

\bibitem{IshizakaStableReduction2026}
Ishizaka, H.: Stable reduction of unresolved dynamics: an operator-theoretic framework. Research note, arXiv:2608.15262 (2026). \url{https://arxiv.org/abs/2608.15262}

\bibitem{LazarZuazua2024}
Lazar, M., Zuazua, E.: Eigenvalue bounds for the Gramian operator of the heat equation. \emph{Automatica} \textbf{164}, 111653 (2024). \url{https://doi.org/10.1016/j.automatica.2024.111653}

\bibitem{LiljegrenSailer2024}
Liljegren-Sailer, B.: Effective error estimation for model reduction with inhomogeneous initial conditions. \emph{Syst. Control Lett.} \textbf{190}, 105840 (2024). \url{https://doi.org/10.1016/j.sysconle.2024.105840}

\bibitem{LionsMagenes1972}
Lions, J.-L., Magenes, E.: \emph{Non-Homogeneous Boundary Value Problems and Applications}, vol.~I. Grundlehren der mathematischen Wissenschaften, vol.~181. Springer, Berlin--Heidelberg (1972). \url{https://doi.org/10.1007/978-3-642-65161-8}

\bibitem{MicchelliRivlin1977}
Micchelli, C.A., Rivlin, T.J.: A survey of optimal recovery. In: Micchelli, C.A., Rivlin, T.J. (eds.) \emph{Optimal Estimation in Approximation Theory}, pp.~1--54. Plenum Press, New York (1977). \url{https://doi.org/10.1007/978-1-4684-2388-4_1}

\bibitem{Mori1965}
Mori, H.: Transport, collective motion, and Brownian motion. \emph{Prog. Theor. Phys.} \textbf{33}, 423--455 (1965). \url{https://doi.org/10.1143/PTP.33.423}

\bibitem{Pazy1983}
Pazy, A.: \emph{Semigroups of Linear Operators and Applications to Partial Differential Equations}. Applied Mathematical Sciences, vol.~44. Springer, New York (1983). \url{https://doi.org/10.1007/978-1-4612-5561-1}

\bibitem{Pruss1993}
Pr\"uss, J.: \emph{Evolutionary Integral Equations and Applications}. Monographs in Mathematics, vol.~87. Birkh\"auser, Basel (1993). \url{https://doi.org/10.1007/978-3-0348-8570-6}

\bibitem{RettbergEtAl2024}
Rettberg, J., Wittwar, D., Buchfink, P., Herkert, R., Fehr, J., Haasdonk, B.: Improved a posteriori error bounds for reduced port-Hamiltonian systems. \emph{Adv. Comput. Math.} \textbf{50}, 100 (2024). \url{https://doi.org/10.1007/s10444-024-10195-8}

\bibitem{SchroederVoigt2023}
Schr\"oder, C., Voigt, M.: Balanced truncation model reduction with a priori error bounds for LTI systems with nonzero initial value. \emph{J. Comput. Appl. Math.} \textbf{420}, 114708 (2023). \url{https://doi.org/10.1016/j.cam.2022.114708}

\bibitem{TraubWasilkowskiWozniakowski1988}
Traub, J.F., Wasilkowski, G.W., Wo\'zniakowski, H.: \emph{Information-Based Complexity}. Academic Press, Boston (1988).

\bibitem{ZhuDominyVenturi2018}
Zhu, Y., Dominy, J.M., Venturi, D.: On the estimation of the Mori--Zwanzig memory integral. \emph{J. Math. Phys.} \textbf{59}, 103501 (2018). \url{https://doi.org/10.1063/1.5003467}

\bibitem{Zwanzig1961}
Zwanzig, R.: Memory effects in irreversible thermodynamics. \emph{Phys. Rev.} \textbf{124}, 983--992 (1961). \url{https://doi.org/10.1103/PhysRev.124.983}

\end{thebibliography}
}
\end{document}